\documentclass[12pt]{article}
\usepackage{amssymb}
\usepackage{amsmath}
\usepackage{latexsym}
\usepackage{enumerate}
\usepackage{hyperref}

\newtheorem{theorem}{Theorem}[section]
\newtheorem{proposition}[theorem]{Proposition}
\newtheorem{corollary}[theorem]{Corollary}
\newtheorem{lemma}[theorem]{Lemma}

\newtheorem{remark}[theorem]{Remark}

\DeclareMathOperator{\id}{id}
\DeclareMathOperator{\vac}{\mathbf{1}}

\DeclareMathOperator{\SL}{SL}
\DeclareMathOperator{\Tr}{Tr}

\DeclareMathOperator{\VB }{VB }
\DeclareMathOperator{\BW}{BW}
\DeclareMathOperator{\Aut}{Aut}

\DeclareMathOperator{\wt}{wt}
\DeclareMathOperator{\Com}{Com}
\newcommand{\qd}{q\frac{\partial}{\partial q}}

\newcommand{\tqed}{\quad $\square$}
\newcommand{\qed}{\quad \square}
\newcommand{\I}{\, i}
\newcommand{\tpi}{2\pi \I}
\newcommand{\td}{\frac{1}{\tpi}\frac{\partial}{\partial \tau}} 

\newcommand{\Trq}[2]{\Tr_{#1}\left ({#2}q^{L(0)-c/24}\right)}

\newcommand{\Zq}[2]{Z_{#1}\left ({#2},q\right)}

\newcommand{\half}{\frac{1}{2}}
\newcommand{\ui}{u_{i}}
\newcommand{\uibar}{\overline{u}_{i}}
\newcommand{\vi}{v_{i}}
\newcommand{\vibar}{\overline{v}_{i}}
\newcommand{\lam}[1]{\lambda^{{#1}}}

\newcommand{\Z}{\mathbb{Z}}
\newcommand{\C}{\mathbb{C}}
\newcommand{\Q}{\mathbb{Q}}
\newcommand{\N}{\mathbb{N}}

\newcommand{\HH}{\mathbb{H}_1}
\newcommand{\omt}{\widetilde{\omega}}

\begin{document}

\title{On Exceptional Vertex Operator (Super) Algebras}
\author{
Michael P. Tuite 
\and 
Hoang Dinh Van\thanks{Supported by Science Foundation Ireland Research Frontiers Programme}
\\
School of Mathematics, Statistics and Applied Mathematics, \\
National University of Ireland Galway \\
University Road, Galway, Ireland}
\maketitle


\maketitle
\abstract{
We consider exceptional vertex operator algebras
and vertex operator superalgebras
with the property that particular Casimir vectors constructed from the primary vectors of lowest conformal weight are Virasoro descendents of the vacuum.  
We show that the genus one partition function and characters for simple ordinary modules must satisfy modular linear differential equations. 
We show the rationality of the central charge and  module lowest weights, modularity of solutions,  the dimension of each graded space is a rational function of the central charge and that the lowest weight primaries generate the algebra.
We also discuss conditions on the reducibility of the lowest weight primary vectors as a module for 
the automorphism group. Finally we analyse  solutions for exceptional vertex operator algebras with primary vectors of lowest weight up to 9 and for 
vertex operator superalgebras with primary vectors of lowest weight up to 17/2. Most solutions can be identified with simple ordinary modules for known algebras but there are also four conjectured algebras generated by weight two primaries and three conjectured extremal vertex operator algebras generated by primaries of weight  3, 4 and 6 respectively. 
}

\section{Introduction}
Vertex Operator Algebras (VOAs) and Super Algebras (VOSAs) have deep connections to Lie algebras, number theory, group theory, combinatorics and Riemann surfaces 
(e.g. \cite{FHL,FLM,Kac1,MN,MT}) and, of course, conformal field theory e.g. \cite{DMS}. 
The classification of VOAs and VOSAs still seems to be a very difficult task, for example,  there is no proof of the uniqueness of the Moonshine module  \cite{FLM}. 
Nevertheless, it would be very useful to be able to characterize 
VOA/VOSAs with interesting properties such as large automorphism groups  (e.g. the Monster group for the Moonshine module), rational  characters, generating vectors etc.  
In \cite{Mat}, Matsuo introduced VOAs of class $\mathcal{S}^n$ with the defining property that the Virasoro vacuum descendents are the only $\Aut(V)$--invariant vectors of weight $k\le n$. 
Thus the the Moonshine module \cite{FLM} is of class $\mathcal{S}^{11}$, the Baby Monster VOA \cite{Ho1} of class $\mathcal{S}^6$ and the level one Kac-Moody VOAs   generated by Deligne's Exceptional Lie algebras $A_1,A_2,G_2,D_4,F_4,E_6,E_7,E_8$ \cite{D} are of class $\mathcal{S}^4$.
\footnote{In fact, the $A_1$ theory is of class  $\mathcal{S}^\infty$ and the $E_8$ theory is of class  $\mathcal{S}^6$.} 

In this paper we consider a refinement and generalization of previous results in \cite{T1,T2} concerning such exceptional VOAs. 
Assuming the VOA is simple and of strong CFT--type (e.g. \cite{MT})
we consider quadratic Casimir vectors $\lambda^{(k)}$ of conformal weight $k=0,1,2,\ldots$ constructed from the primary vectors of lowest conformal weight $l\in \N$.
We say that a VOA is Exceptional of Lowest Primary Weight $l$ if $\lambda^{(2l+2)}$ is a Virasoro vacuum descendent. 
Every VOA of class $\mathcal{S}^{2l+2}$ with lowest primary weight $l$  is Exceptional  but the converse is not known to be true.
We show, using Zhu's  theory for genus one correlation functions \cite{Z}, that for an Exceptional VOA of lowest primary weight $l$, the partition function and the characters for simple ordinary VOA  modules  satisfy a Modular Linear Differential Equation (MLDE) of order at most $ l+1$. 
Given that order of the MLDE  is exactly $l+1$ (which is verified for all $l\le 9$)
we show that the central charge $c$ and module lowest weights $h$ are rational, the MLDE solution space is modular invariant and the dimension of each VOA graded space is a rational function of $c$.
Subject to a further indicial root condition (again verified for all $l\le 9$) we  show that an Exceptional VOA is generated by its primary vectors of lowest weight $l$.

We also consider other properties that arise from genus zero correlation functions for all $l$. 
Assuming the VOA is of class $\mathcal{S}^{2l+2}$ this leads to  conditions on the reducibility of the lowest weight $l$ primary space as a module for the VOA automorphism group.

A similar analysis is carried out for Exceptional VOSAs of lowest primary weight $l\in \N+\half$ for which $\lambda^{(2l+1)}$ is a Virasoro vacuum descendent. 
Using a twisted version of Zhu theory \cite{MTZ} we obtain a Twisted MLDE of order at most $l+\half $ which is satisfied by the partition function and simple ordinary VOA  module characters. 
This differential equation leads to a similar set of general results to those for VOAs. Likewise, we can consider genus zero correlation functions for all  $l\in \N+\half$ leading to conditions on the reducibility of the space of the space of weight $l$ primaries as a module for the VOSA automorphism group. 

The paper also summarizes  rational $c,h$ solutions to the MLDE for all $l\le 9$ and the Twisted MLDE for all $l\le \frac{17}{2}$. In most cases we can identify a VOA/VOSA with the requisite properties. These include a number of  special VOA/VOSA constructions, some commutant VOSA constructions,  some Virasoro minimal model simple current extensions and ${\mathcal W}$--algebras. We also present evidence for four candidate/conjectured  VOAs with simple Griess algebras for $l=2$ and three extremal VOAs for $l=3,4,6$. All the VOSA solutions found can be identified with known theories.

\section{Vertex Operator (Super) Algebras}
\label{vos}
We review some aspects of Vertex Operator Super Algebra  theory (e.g. 
\cite{FHL,FLM,Kac1,MN,MT}).
A Vertex Operator Superalgebra (VOSA) is a quadruple 
$(V,Y(\cdot,\cdot),\vac,\omega)$ with a $\Z_2$-graded vector space $V=V_{\bar{0}}\oplus V_{\bar{1}}$
with parity $p(u)=0$ or $1$ for $u\in V_{\bar{0}}$ or $V_{\bar{1}}$ respectively. $(V,Y(\cdot,\cdot),\vac,\omega )$ is called a Vertex Operator Algebra (VOA) when $V_{\bar{1}}=0$. 
\medskip

$V$ also has a $\half\Z$-grading with $V=\bigoplus _{r \in \half\Z}V_{r}$  
with $\dim V_{r}<\infty$.
   $\vac\in V_{0}$ is the vacuum vector and 
$\omega\in V_{2}$ is called the conformal vector.
$Y$ is a linear map $Y:V\rightarrow \mathrm{End}(V)[[z,z^{-1}]]$ for formal
variable $z$ giving a vertex operator  
\begin{equation}
Y(u,z)=\sum_{n\in \Z}u(n)z^{-n-1},  
\label{Ydefn}
\end{equation}
 for every $u\in V$.
The linear operators (modes) $u(n):V\rightarrow V$ satisfy creativity 
\begin{equation}
Y(u,z)\vac = u +O(z),
\label{create}
\end{equation}
and lower truncation 
\begin{equation}
u(n)v=0,
\label{lowertrun}
\end{equation}
for each $u,v\in V$ and $ n\gg 0$.
For the conformal vector $\omega$ 
\begin{equation}
Y(\omega ,z)=\sum_{n\in \Z}L(n)z^{-n-2},  \label{Yomega}
\end{equation}
where $L(n)$ satisfies the Virasoro algebra for some central charge $c$ 
\begin{equation}
[ L(m),L(n)]=(m-n)L(m+n)+\frac{c}{12}(m^{3}-m)\delta_{m,-n}\id_V.
\label{Virasoro}
\end{equation}
Each vertex operator satisfies the translation property 
\begin{equation}
Y(L(-1)u,z)=\partial_z Y(u,z).  
\label{YL(-1)}
\end{equation}
The Virasoro operator $L(0)$ provides the $\half\Z$-grading with $L(0)u=\wt(u)u$ for 
$u\in V_{r}$ and with weight $\wt(u)=r\in \Z+\half p(u)$. 
Finally, the vertex operators satisfy the Jacobi identity
\begin{eqnarray*}
\notag
& z_0^{-1}\delta\left( \frac{z_1 - z_2}{z_0}\right) Y (u, z_1 )Y(v , z_2)    
  -(-1)^{p(u)p(v)}z_0^{-1} \delta\left( \frac{z_2 - z_1}{-z_0}\right) Y(v, z_2) Y(u , z_1 ) 
&\\
\label{VOAJac}
& 
= z_2^{-1}    
\delta\left( \frac{z_1 - z_0}{z_2}\right)
Y \left( Y(u, z_0)v, z_2\right).  &
\end{eqnarray*} 
with $\delta\left(\frac{x}{y}\right)=\sum_{r\in \Z} x^{r}y^{-r}$.

These axioms imply $u(n)V_r\subset V_{r-n+wt(u)-1}$ for $u$ of weight $\wt(u)$. 
They also imply
locality, skew-symmetry, associativity and commutativity:
\begin{eqnarray}
&&(z_{1}-z_{2})^N
Y(u,z_{1})Y(v,z_{2}) 
= (-1)^{p(u)p(v)}(z_{1}-z_{2})^N
Y(v,z_{2})Y(u,z_{1}),
\notag
\\
\label{Local}\\
&&Y(u,z)v = (-1)^{p(u)p(v)} e^{zL(-1)}Y(v,-z)u,
\label{skew}
\\
&&(z_{0}+z_{2})^N Y(u,z_{0}+z_{2})Y(v,z_{2})w = (z_{0}+z_{2})^N Y(Y(u,z_{0})v,z_{2})w,
\label{Assoc}\\
&&u(k)Y(v,z)-(-1)^{p(u)p(v)} Y(v,z)u(k)
=\sum_{j\ge 0}\binom{k}{j}
Y(u(j)v,z)z^{k-j},\notag
\\ \label{Comm}
\end{eqnarray}
for $u,v,w\in V$ and integers $N\gg 0$  \cite{FHL},  \cite{Kac1},  \cite{MT}. 

We define an invariant symmetric bilinear form $\langle\ ,\ \rangle $ on $V$ by
\begin{equation}
\left \langle Y\left(e^{zL(1)}\left(-z^{-2}\right)^{L(0)}w, z^{-1}\right)u,v\right \rangle
=(-1)^{p(u)p(w)}\langle v,Y(w,z)v\rangle.  \label{Yinvariant}
\end{equation}%
for all $u,v,w\in V$  \cite{FHL}.
$V$ is said to be of \emph{CFT--type} if $V_0=\C \vac$ and of \emph{Strong CFT-Type} if additionally $L(1)V_1 =0$ in which case $\langle\ ,\ \rangle$, with normalization
$\langle \mathbf{1},\mathbf{1}\rangle =1$, is unique \cite{Li}.
Furthermore, $\langle\ ,\ \rangle$ is invertible if $V$ is simple. 
All VOSAs in this paper are assumed to be of this type.

Every VOSA contains the subVOA $V_{\omega}$ generated by the Virasoro vector $\omega$ with Fock basis
of vacuum descendents of the form
\begin{eqnarray}
L(-n_1)L(-n_2)\ldots L(-n_k)\vac,
\label{eq:VirFock}
\end{eqnarray}
for $n_i\ge 2$. 
$\langle\ ,\ \rangle $ is singular on $\left(V_{\omega}\right)_n$ iff the central charge  is 
\begin{equation}
c_{p,q}=1-6\frac{(p-q)^{2}}{pq},  \label{cpq}
\end{equation}
 for coprime integers $p,q\geq 2$ and $n=(p-1)(q-1)$ \cite{Wa}.
The Virasoro minimal model VOA $L(c_{p,q},0)$ is the quotient of $V_{\omega}$ by the radical of $\langle\ ,\ \rangle $.
$L(c_{p,q},0)$ has a finite number of simple ordinary $V$--modules $L(c_{p,q},h_{r,s})\cong L(c_{p,q},h_{q-r,p-s})$ (e.g. \cite{DMS})  with lowest weight 
\begin{eqnarray}
h_{r,s}={\frac { \left( pr-qs \right) ^{2}- \left( p-q \right) ^{2}}{4pq}},
\label{eq:hrs}
\end{eqnarray}
for $r=1,\ldots,q-1$ and $s=1,\ldots ,p-1$.

\section{Quadratic Casimirs and Genus One Zhu Theory}
\subsection{Quadratic Casimirs}
\label{sect:Cas}
Let $(V,Y(\cdot,\cdot),\vac,\omega)$ be a simple VOA of strong CFT--type with unique invertible bilinear form $\langle\ ,\ \rangle$. 
Let $\Pi_l$ denote the space of  primary vectors of lowest weight $l\ge 1$  i.e. $L(n)u=0$ for all $n>0$ for $u\in \Pi_l$.
Choose a $\Pi_l$--basis $\{\ui\}$  for $i=1,\ldots ,p_l=\dim \Pi_l$   with dual basis $\{\uibar\}$  i.e. $\langle\ui ,\overline{u}_{j} \rangle=\delta_{ij}$.  Define quadratic Casimir vectors $\lam{(n)}$ for $n\ge 0$ by \cite{Mat,T1,T2}
\begin{eqnarray}
\lam{(n)}=\ui(2l-n-1)\uibar\in V_{n},
\label{eq:lamdef}
\end{eqnarray}
where the repeated $i$ index is summed from 1 to $p_l$ (here and below).
Then 
\[
\lam{(0)}=\ui(2l-1)\uibar=\left(-1\right)^l\langle \ui,\uibar\rangle\vac=\left(-1\right)^lp_l \vac.
\]
Furthermore, if $l>1$ then $\dim V_{1}=0$ and hence $\lam{(1)}=0$ whereas for $l=1$  the Jacobi identity implies 
$
\lam{(1)}=\ui(0)\uibar=-\uibar(0)\ui=0
$ \cite{T1}.
Thus we find 
\begin{lemma}\label{lemma:lambda01}
$\lam{(0)}=\left(-1\right)^lp_l \vac$ and $\lam{(1)}=0$.
\end {lemma}
Since the $\Pi_l$ elements are primary then for all $m>0$
\begin{equation}
L(m)\lambda^{(n)}=\left(n-m+l(m-1)\right)\lambda^{(n-m)}.
\label{eq:LMlambda}
\end{equation}
Suppose that $\lambda^{(n)}\in  V_{\omega}$ then \eqref{eq:LMlambda} implies that $\lambda^{(m)}\in V_{\omega}$ 
for all $m\le n$. Furthermore, since $\langle\ ,\ \rangle$ is invertible 
we have \cite{Mat}
\begin{lemma}\label{lemma:lambda}
If $\lambda^{(n)}\in V_{\omega}$ then $\lambda^{(n)}$ is uniquely determined. 
\end {lemma}
Thus if $\lambda^{(2)}\in V_{\omega}$ then  $\lambda^{(2)}=\kappa L(-2)\vac$ for
 some $\kappa$ so that 
$\langle L(-2)\vac ,\lambda^{(2)}\rangle =\kappa\frac{c}{2}$. 
But \eqref{Yinvariant} and \eqref{eq:LMlambda} imply
$\langle L(-2)\vac,\lambda ^{(2)}\rangle =
\langle \vac,L({2})\lambda ^{(2)}\rangle =\left(-1\right)^l p_l$ so that for $c\neq c_{2,3}=0$ (c.f. \eqref{cpq})
\begin{equation}
\lambda^{(2)}= p_l \frac{2 \left(-1\right)^l l}{c}L(-2)\vac.
\label{eq:lambda2}
\end{equation}
Similarly, if  $\lambda^{(4)}\in V_{\omega}$ and $c\neq 0,c_{2,5}=-22/5$ then  \cite{Mat, T1,T2}
\begin{equation}
\lam{(4)}=
p_l\frac{2 \left( -1 \right) ^{l}l\left( 5 l+1 \right)}{c \left( 5c + 22 \right) } L(-2)^2\vac
+ p_l\frac{3 \left( -1 \right) ^{l}l\left(  c-2 l+4 \right)}{c \left( 5c + 22 \right) } L(-4)\vac.
\\  \label{eq:lambda4}
\end{equation}
These examples illustrate a general observation:
\begin{lemma}
\label{lem:ratlam}
Each coefficient in the expansion of $\lam{(n)}\in  V_{\omega}$ in a basis of Virasoro Fock vectors is of the form $p_l \,r(c)$ for $r(c)$ a rational  function of $c$.
\end{lemma}

\subsection{Genus One Constraints from Quadratic \\ Casimirs}
Define genus one partition and 1-point correlation functions for $u\in V$ by 
\begin{eqnarray}
Z_V(q)&=&\Trq{V}{}=q^{-c/24}\sum_{n\geq 0}\dim V_n \,q^n,\label{eq:Z(q)}
\\
Z_{V}(u,q)&=&\Trq{V}{o(u)},\label{eq:1pt}
\end{eqnarray}
for formal parameter $q$ and for \lq zero mode\rq\ $o(u)=u({\wt{(u)}-1}):V_n\rightarrow V_n$ for homogeneous $u$. 
By replacing $V$ by a simple ordinary $V$--module $N$ (on which $L(0)$ acts semi-simply e.g. \cite{FHL, MT}) these definitions may be extended to $N$ graded characters $Z_{N}(q)$ and 1-point functions $Z_{N}(u,q)$. Thus 
\begin{equation}
Z_N(q)=\Trq{N}{}=q^{h-c/24}\sum_{n\geq 0}\dim N_n \,q^n,\label{eq:ZN(q)}
\end{equation}
where $h$ denotes the lowest weight of $N$. Zhu also  introduced an isomorphic VOA $(V, Y[\cdot,\cdot], \vac,\omt)$ with \lq square bracket\rq\  vertex operators   
\begin{equation}
Y[u,z]\equiv Y\left(e^{zL(0)}u,e^{z}-1\right)
=\sum\limits_{n\in \Z}u[n]z^{-n-1},
\label{eq:square}
\end{equation}%
for  Virasoro vector $\omt=\omega -c/24\vac$  with modes $\{L[n]\}$. 
$L[0]$ defines an alternative $\Z$ grading with $V=\bigoplus_{k\geq 0}V_{[k]}$
where $L[0]v=\wt[v]v$ for $\wt[v]=k$ for $v\in V_{[k]}$.  
Zhu obtained a reduction formula for the 2-point correlation function $\Zq{N}{Y[u,z]v}$ for $u,v\in V$
in terms of the elliptic Weierstrass function
\begin{equation}
P_{m}(z)=\frac{1}{z^{m}}+(-1)^{m}\sum_{n\ge m}\binom{n-1}{m-1}E_{n}(q)\, z^{n-m},
\label{eq:Pm}
\end{equation}
 for $m\ge 1$ and with Eisenstein series $E_{n}(q)=0$ for odd $n$  and 
\begin{equation}
E_{n}(q)=-\frac{B_{n}}{n!}+\frac{2}{(n-1)!}
\sum_{k\geq 1}\frac{k^{n-1}q^{k}}{1-q^{k}}, \label{eq:En}
\end{equation}
for even $n$  with $B_{n}$  the $n$th Bernoulli number. $P_{m}(z)$ converges absolutely and uniformly on compact subsets of the domain
$|q| < |e^{z} | < 1$. 
$E_n(q)$ is a modular form of weight $n$ for $n\ge 4$ and  $E_2(q)$ is a quasi-modular form of weight 2 i.e. letting $q=\exp(\tpi \tau)$ for $\tau\in \HH$
\begin{equation}
E_n\left(\frac{\alpha\tau +\beta}{\gamma\tau +\delta}\right)=\left( \gamma\tau +\delta\right)^{n}E_n(\tau)-\frac{\gamma(\gamma\tau +\delta)}{\tpi}\delta_{n2},
\label{eq:Ekmod}
\end{equation}
for $\bigl(\begin{smallmatrix}
\alpha&\beta\\ \gamma&\delta
\end{smallmatrix} \bigr)\in \SL(2,\Z)$ \cite{Se}. 
 We then have \cite{Z}
\begin{proposition}[Zhu]
\label{prop:Zhu}
Let $N$ be a simple ordinary $V$--module.
\begin{eqnarray*}
\Zq{N}{ Y[u,z]v}
&=&
\Trq{N}{o(u)o(v)}\\
&&+\sum_{m\ge 0} P_{m+1}(z) \Zq{N}{u[m]v}.
\end{eqnarray*}
\end{proposition}
Taking $u=\omt$ and noting that $o(\omt)=L(0)-c/24$ we obtain:
\begin{corollary}
\label{thm:ZhuLk}
The 1-point function of a Virasoro descendent $L[-k]v$ is
\begin{eqnarray*} 
\Zq{N}{ L[-2]v} &=&\left( \qd +\wt[v]E_2(q)\right)\Zq{N}{v}\\
&& + \sum_{s\ge 1} E_{2s+2}(q) \Zq{N}{L[2s]v},\\
\Zq{N}{ L[-k]v} &=& (-1)^{k} \sum_{r\ge 0}\binom{k+r-1}{k-2} E_{k+r}(q) \Zq{N}{L[r]v},
\end{eqnarray*}
for all  $ k\ge 3$.
\end{corollary}
\medskip

Let us now consider a simple VOA $V$ of strong CFT-type  
with lowest weight $l\ge 1$ Virasoro primary vectors $\Pi_l$ so that
\begin{eqnarray}
Z_{V}(q)=Z_{V_{\omega}}(q)+O\left(q^{l+c/24}\right).
\label{eq:Zlead}
\end{eqnarray}
Let $\{\ui\}$ and $\{\uibar\}$ be a basis and dual basis  for $\Pi_l$. 
Apply Proposition~\ref{prop:Zhu} to 
\begin{eqnarray}
\Zq{N}{ Y[\ui,z]\uibar}=\sum_{n\ge 0}\Zq{N}{\lam{[n]}}z^{n-2l},
\label{Z2pt}
\end{eqnarray}
(for Casimir vector $\lam{[n]}\in V_{[n]}$ in square bracket modes) to find 
\begin{eqnarray}
\sum_{n\ge 0}\Zq{N}{\lam{[n]}}z^{n-2l}&=& 
\Trq{N}{o(\ui)o(\uibar)}
\notag\\
&&
+\sum_{m=0}^{2l-1} P_{m+1}(z) \Zq{N}{\lam{[2l-m-1]}}
.\quad
\label{eq:lamrec}
\end{eqnarray}
Equating $z^{n-2l}$ coefficients results in recursive identities for $\Zq{N}{\lam{[n]}}$ for $n\ge 2l$.
In particular, equating the  $z^2$ coefficients implies
\begin{proposition}
\label{prop:lamrec}
$\Zq{N}{\lam{[2l+2]}}$ satisfies the recursive identity
\begin{eqnarray}
\Zq{N}{\lam{[2l+2]}}&=&
\sum_{r=0}^{l-1} 
\binom{2l-2k+1}{2}
E_{2l-2k+2}(q)
\Zq{N}{\lam{[2k]}}.\quad 
\label{Zrec2}
\end{eqnarray}
\end{proposition}

\section{Exceptional VOAs}
\label{sec:EVOAs}
Consider a simple VOA of strong CFT--type  with primary vectors  of lowest weight $l\ge 1$ for which $\lam{(2l+2)}\in V_{\omega}$ (or equivalently, $\lam{[2l+2]}\in  V_{\omt}$). 
We also assume that $ \left(V_{\omega}\right)_{2l+2}$ contains no Virasoro singular vector i.e. $c\neq c_{p,q}$ for $(p-1)(q-1)\le 2l+2$. 
We call such a VOA  an \emph{Exceptional VOA of Lowest Primary Weight $l$}. 
\eqref{eq:LMlambda} implies  $\lam{(2k)}\in V_{\omega}$ (and $\lam{[2k]}\in V_{\omt}$) for all $k\le l$. 
\begin{proposition}  
\label{prop:lam2k}
Let $\lam{[2k]}\in V_{\omt}$. Then for a simple ordinary $V$--module $N$
\begin{eqnarray}
\Zq{N}{\lam{[2k]}}=
\sum_{m=0}^{k} f_{k-m}(q,c)D^{m}Z_N(q),
\label{eq:lamD}
\end{eqnarray}
where $D$ is the Serre modular derivative defined for $m\ge 0$ by
\begin{equation}
D^{m+1} Z_N(q)=\left(\qd +2mE_2(q)\right)D^{m}Z_N(q).
\label{eq:Dm}
\end{equation} 
 $f_{m}(q,c)$ is a modular form of weight $2m$ whose coefficients over the ring of Eisenstein series are of the form $p_l \,r(c)$ for  a rational function $r(c)$.
\end{proposition}

\noindent \textbf{Proof.}
\eqref{eq:lamD} follows from Corollary~\ref{thm:ZhuLk} by induction in the number of Virasoro modes where the $D^kZ_N(q)$  term  arises from a $L[-2]^k\vac$ component in $\lam{[2k]}$.
The coefficients of $f_{m}(q,c)$ over the ring of Eisenstein series are of the form $p_l\, r(c)$ for a rational function $r(c)$ from Lemma~\ref{lem:ratlam}.\tqed
\
Applying Proposition~\ref{prop:lam2k} to the the recursive identity \eqref{Zrec2} implies $Z_N(q)$ satisfies a Modular Linear Differential Equation (MLDE) \cite{Mas1} 
\begin{proposition}
\label{prop:MLDE}
Let $V$ be an Exceptional VOA of lowest primary weight $l$.  $Z_N(q)$ for each  simple ordinary $V$--module $N$  satisfies a MLDE of order $\le l+1$
\begin{equation}
\sum_{m=0}^{l+1} g_{l+1-m}(q,c)D^{m}Z(q)=0,
\label{eq:MLDE}
\end{equation}
where  
 $g_{m}(q,c)$ is a modular form of weight $2m$ whose coefficients over the ring of Eisenstein series are rational functions of $c$.
\end{proposition}
$g_{0}(q,c)=g_0(c)$ is independent  of $q$ since it is a modular form of weight $0$.
For $g_{0}(c)\neq 0$, the  MLDE  \eqref{eq:MLDE} is of order $l+1$ with a regular singular point at $q=0$ so that  Frobenius-Fuchs theory  concerning the $l+1$ dimensional solution space  ${\mathcal F}$ applies e.g. \cite{Hi,I}. 
Any solution $Z(q)\in {\mathcal F}$ is holomorphic in $q$ for $0<|q|<1$
since the MLDE coefficients $g_{m}(q,c)$ are holomorphic for $|q|<1$.
We may thus view each solution as a function of $\tau\in \HH$ for  $q=e^{\tpi \tau}$. 

 Using the quasi-modularity of $E_2(\tau)$ and  \eqref{eq:Dm} with $\qd=\td$, it follows that
 for all 
$\bigl(\begin{smallmatrix}
\alpha&\beta\\ \gamma&\delta
\end{smallmatrix} \bigr)\in \SL(2,\Z)$,
$Z\left(\frac{\alpha\tau +\beta}{\gamma\tau +\delta}\right)$ is also a  solution of the MLDE since $g_{l+1-m}(q,c)$ is a modular form of weight $2l+2-2m$. 
Thus $T:\tau\rightarrow \tau+1$ has a natural action on ${\mathcal F}={\mathcal F}_1\oplus\ldots \oplus {\mathcal F}_r$ for $T$ eigenspaces ${\mathcal F}_i$ with monodromy eigenvalue $e^{\tpi x}$. $x$ is a root of the  indicial equation
\begin{eqnarray}
\sum_{m=0}^{l+1} g_{l+1-m}(0,c)\prod_{s=0}^{m-1}(x-\frac{1}{6}s) =0.
\label{eq:indicial}
\end{eqnarray}
If  $x_1=x_2\mod \Z$, for roots $x_1,x_2$, they determine the same monodromy eigenvalue. Let $x_i$ denote the  indicial root with least real part for a given monodromy eigenvalue.
Then ${\mathcal F}_i$ has a basis of the form
\begin{equation}
f_{i}^n(\tau)=\phi_{i}^1(q)+\tau \phi_{i}^2(q)+\ldots + \tau^{n-1} \phi_{i}^{n}(q),
\label{eq:logbasis}
\end{equation}
for $q$--series (which are holomorphic on $0<|q|<1$) of the form
\[
\phi_{i}^n(q)=q^{x_i}\sum_{k\ge 0} a_{ik}^n q^k,\quad 1\le n\le \dim {\mathcal F}_i.
\]
Logarithmic solutions (with  a $\tau^j$ factor for $j>0$) occur if the same root occurs multiple times or, possibly, if two roots differ by an integer.
However, every graded character $Z_{N}(q)$ for a simple  ordinary module with  lowest weight $h$ has a pure $q$--series with indicial root $x=h-c/24$ from \eqref{eq:ZN(q)}. 

We now sketch a proof that the central charge $c$ is rational following \cite{AM} (which is extended to logarithmic solutions \eqref{eq:logbasis} in \cite{Miy}). 
Suppose $c\notin \Q$ and consider $\phi\in \Aut(\C)$ such that $\tilde{c}=\phi(c)\neq c$.
Then $Z_V(\tau,\tilde{c})$ is a solution to the MLDE \eqref{eq:MLDE} found by replacing $c$ by $\tilde{c}$.
But since the coefficients in the $q$--expansion of $Z_V(\tau,c)$ are integral we have 
\[
Z_V(\tau,\tilde{c})=q^{(\tilde{c}-c)/24} Z_V(\tau,c).
\]
%
Applying the  modular transformation $S:\tau\rightarrow -1/\tau$ we find
\begin{equation}
Z_V\left(-\frac{1}{\tau},\tilde{c}\right)=\exp\left(-\frac{\pi\I(\tilde{c}-c)}{12\tau}\right) Z_V\left(-\frac{1}{\tau},c\right).
\label{eq:ZVS}
\end{equation}
But  $Z_V(-1/\tau,c)$ satisfies \eqref{eq:MLDE} and $Z_V(-1/\tau,\tilde{c})$ satisfies \eqref{eq:MLDE} with $c$ replaced by $\tilde{c}$  and thus both are of the form \eqref{eq:logbasis}. 
Analysing  \eqref{eq:ZVS} along rays $\tau=re^{\I \theta}$ in the limit $r\rightarrow \infty$ with $0 <\theta< \pi$ a contradiction results unless  $\tilde{c}=c$. Hence $c\in \Q$ \cite{AM,Miy}.  
Similarly, the lowest conformal weight $h$ of a simple ordinary module $N$ is rational. Altogether we have
\begin{proposition}
\label{prop:Zsol}
Let $V$ be an Exceptional VOA of lowest primary weight $l\ge 1$  and central charge $c$ and let $N$ be a simple ordinary $V$--module of lowest weight $h$. Assuming $g_{0}(c)\neq 0$ in the MLDE  \eqref{eq:MLDE} then
\begin{enumerate}
\item[(i)]   $Z_{N}(q)$ is  holomorphic for $0<|q|<1$. 
\item[(ii)]   
$Z_{N}\left(\frac{\alpha\tau +\beta}{\gamma\tau +\delta}\right)$ is a solution of the MLDE  for all
 $\bigl(\begin{smallmatrix}
\alpha&\beta\\ \gamma&\delta
\end{smallmatrix} \bigr)\in \SL(2,\Z)$ 
viewed as a function of $\tau\in \HH$ for $q=e^{\tpi\tau}$. 
\item[(iii)] The central charge $c$  and the lowest conformal weight $h$ are rational. 
\end{enumerate}
\end{proposition}

Consider the general solution  with indicial root $x=c/24$ of the form $Z(q)=q^{-c/24}\sum_{n\ge 0} a_n q^n$.
Substituting into the MLDE we 
obtain a linear equation in  $a_0,\ldots ,a_n$ for each $n$. 
This can be  iteratively solved for $a_n$  provided the coefficient of $a_n$ is non-zero. 
This coefficient may vanish if $x=m-c/24$ is an indicial root for some integer $m> 0$. Hence we have
\begin{proposition}
\label{prop:Zgen}
Let $V$ be an Exceptional VOA of lowest primary weight $l\ge 1$  and central charge $c$. 
Suppose $g_{0}(c)\neq 0$ and that $m<l$ for any indicial root of the form $x=m-c/24$. Then
\begin{enumerate}
\item[(i)]   $Z_V(q)$ is the unique $q$--series solution of the MLDE  obeying
\eqref{eq:Zlead}.
\item[(ii)] $\dim V_n$ is a rational function of $c$ for each $n\ge 0$.
\item[(iii)]  $V$  is generated by the space of  lowest weight primary vectors $\Pi_l$.

\end{enumerate}
\end{proposition}

\noindent \textbf{Proof.}
(i) The $x=-c/24$ solution  $Z(q)=q^{-c/24}\sum_{n\ge 0} a_n q^n$ is  determined by $a_0$ and $a_m$ for any indicial root(s) of the form $x=m-c/24$ for $m>0$. Thus the partition function is uniquely determined by the $l$ Virasoro leading terms \eqref{eq:Zlead} under the assumption that $m<l$. 

(ii) The modular forms $g_m(q,c)$ of the MLDE of Proposition~\ref{prop:MLDE} have $q$--expansions whose coefficients are rational functions of $c$. Hence solving iteratively  it follows that 
$a_n=\dim V_n$ is a rational function of $c$. 

(iii) Let $V_{\langle \Pi_l\rangle}\subseteq V$ be the subalgebra generated by the lowest weight  primary vectors $\Pi_l$.  
But $\omega\in V_{\langle \Pi_l\rangle}$ from  \eqref{eq:lambda2} so that $ V_{\langle \Pi_l\rangle}$ is a VOA of central charge $c$. 
Furthermore, since $\lam{(2l+2)}\in V_{\langle \Pi_l\rangle}$, the subVOA is an Exceptional VOA of lowest primary weight $l$. 
Hence $Z_{V_{\langle \Pi_l\rangle}}(q)$ obeys the same MLDE as $Z_{V}(q)$.
From (i) it follows that $
Z_{V_{\langle \Pi_l\rangle}}(q)=Z_{V}(q)$ implying $V_{\langle \Pi_l\rangle}= V$.
\tqed
\
\begin{remark} \label{rem_g0}
Note that $g_{0}(c)\neq 0$ provided $\lam{(2l+2)}$ contains an $L(-2)^{l+1}\vac$ component. We conjecture that such a component exists for all $l$. 
We further conjecture that $m<l$ for any indicial root of the form $x=m-c/24$ for all $l$. These properties are  verified for all $l\le 9$ in  Section~\ref{sect: VOAex}.
\end{remark}

\subsection{ Exceptional VOAs with $p_l=1$}
\label{sect:pl1}
Let $V$ be a simple VOA of strong CFT type generated by one primary vector $u$ of lowest  weight $l$ with  dual $\overline{u}=u/\langle u,u\rangle$. 
Consider the commutator \eqref{Comm}
\begin{eqnarray}
\left[ u(m),Y(u,z)\right] &=& \sum_{j\ge 0}\binom{m}{j}Y\left (u(j)u,z \right)z^{m-j}
\notag
\\
&=& \langle u,u\rangle \sum_{k= 0}^{2l-1}\binom{m}{2l-k-1}Y\left (\lambda^{(k)},z\right)z^{m+k+1-2l},\quad
\label{uYcom}
\end{eqnarray}
using \eqref{eq:lamdef}. 
Suppose that  $\lam{(2l-1)}\in V_{\omega}$ so that $\lambda^{(k)}\in V_{\omega}$ for $0\le k\le 2l-1$ which implies the RHS of \eqref{uYcom} is expressed in terms of  Virasoro modes.
Thus \eqref{uYcom} defines a ${\mathcal W}(l)$  algebra VOA generated by $u$ e.g. \cite{BS}. 
The further  condition $\lam{(2l+2)}\in V_{\omega}$ constrains $c$ to specific rational values.
 
We consider two infinite families of Exceptional ${\mathcal W}(l)$--VOAs. One is of $AD$-type, from the $ADE$ series of \cite{CIZ}, given by the simple current extension of a minimal model $L\left(c_{p,q},0\right)$  by an irreducible module $L\left(c_{p,q},l \right)$ with
\begin{equation}
l=h_{1,p-1}=\frac{1}{4}(p-2)(q-2)\in\N,
\label{eq:hcur}
\end{equation}
for $h_{r,s}$ of \eqref{eq:hrs} i.e. for any coprime pair $p,q$ such that $p$ or $q= 2 \mod 4$. Then \eqref{uYcom} is compatible with the Virasoro fusion rule (e.g. \cite{DMS})
\[
L\left(c_{p,q},h_{1,p-1}\right)\times L\left(c_{p,q},h_{1,p-1}\right) = L\left(c_{p,q},0\right).
\]
Furthermore, since
\[
2l+2=(p-1)(q-1)-\frac{1}{2}(pq-6)<(p-1)(q-1),
\]
it follows that $\left(V_{\omega}\right)_{2l+2}$ contains no Virasoro singular vectors.
Hence
\begin{proposition}
\label{prop:ADVOA}
For a minimal model with $h_{1,p-1}\in\N$ there exists an Exceptional VOA with one primary vector of  lowest  weight $l=h_{1,p-1}$ of $AD$-type
\begin{equation}
V=L\left(c_{p,q},0\right)\oplus L\left(c_{p,q},h_{1,p-1}\right).
\label{eq:Vsimcur}
\end{equation}
\end{proposition}

A second infinite family of ${\mathcal W}(l)$--VOAs for $l=3k$ for $k\ge 1$ is given in \cite{BFKNRV, F}.  
A more complete VOA description of this construction will appear elsewhere \cite{T3}. 
${\mathcal W}(3k)$ is of central charge $c_k=1-24k$ and contains a unique 
Virasoro primary vector of weight $h_n=(n^2-1)k$ for each $n\ge 1$. 
The corresponding Virasoro Verma module contains a unique singular vector of weight $h_n+n^2$ so that the 
partition function is \cite{F}
\begin{eqnarray}
Z_{{\mathcal W}(3k)}(q)&=&\sum_{n\ge 1}\frac{q^{-c_k/24}}{\prod_{m\ge 0}(1-q^m)}\left (q^{h_{n}} -q^{h_{n}+n^2}\right)\notag
\\
&=&\frac{1}{2\eta(q)}\sum_{n\in\Z}\left (q^{n^2 k}-q^{n^2(k+1)}\right).
\label{eq:ZW3k}
\end{eqnarray}
This VOA is generated by the lowest weight primary of weight $l=h_2=3k$ 
$\lam{(2l+2)}\in \left(V_{\omega}\right)_{2l+2}$ requires that
$h_3=8k>2l+2$ i.e. $k>1$. Thus we find
\begin{proposition}
\label{prop:W3kVOA}
For each $k\ge 2$ there exists an Exceptional VOA ${\mathcal W}(3k)$ with one primary vector of  lowest  weight $3k$ and central charge $c_k=1-24k$.
\end{proposition}
\begin{remark} \label{Walg}
We conjecture that the two VOA series  of Propositions~\ref{prop:ADVOA} and \ref{prop:W3kVOA} are the only Exceptional VOAs for which $p_l=1$. 
\end{remark}

\section{Genus Zero Constraints from Quadratic Casimirs}
\label{sect:genus 0}  
We next consider how an Exceptional VOA is also subject to local genus zero constraints following an approach originally described for $l=1,2$ in \cite{T1,T2}. 
Let $V$ be a simple VOA of strong CFT-type with lowest primary weight $l\ge 1$.  Let $\Pi_l$ be the vector space of $p_l$ primary vectors of weight $l$ with basis $\{\ui\}$ and dual basis $\{ \uibar\}$. 
Define the genus zero correlation function 
\begin{eqnarray}
F(a,b;x,y)=\langle a,Y(\ui,x)Y(\uibar,y)b\rangle,\label{eq:Fab}
\end{eqnarray}
for $a,b\in \Pi_l$. Note that $F(a,b;x,y)$ is constructed locally from $\Pi_l$ alone.
Locality \eqref{Local}, associativity \eqref{Assoc} and lower truncation \eqref{lowertrun} give
\begin{proposition}\label{prop:FabG}
$F(a,b;x,y)$ is determined  by a rational function
\begin{equation}
F(a,b;x,y)=\frac{G(a,b;x,y)}{x^{2l}y^{2l}(x-y)^{2l}},
\label{FabG}
\end{equation}
for $G(a,b;x,y)$ a symmetric homogeneous polynomial in $x,y$
 of degree $4l$.
\end{proposition}
$F(a,b;x,y)$ can be considered as a rational function on the genus zero Riemann sphere and expanded in a various domains to obtain the $2l+1$ independent parameters determining $G(a,b;x,y)=\sum_{r=0}^{4l}A_{r}x^{4l-r}y^{r}$ where $A_{r}=A_{4l-r}$. In particular, 
we expand in $\xi=-y/(x-y)$ using skew-symmetry \eqref{skew}, translation \eqref{YL(-1)} and invariance of $\langle\ ,\ \rangle$ to find 
\begin{eqnarray}
y^{2l} F(a,b;x,y) &=&y^{2l}\langle a,Y(\ui,x)e^{yL_{-1}}Y(b,-y)\uibar\rangle
\notag \\
&=&y^{2l}\langle a,e^{yL_{-1}}Y(\ui,x-y)Y(b,-y)\uibar\rangle  \notag
\\
&=&y^{2l}\langle a,Y(\ui,x-y)Y(b,-y)\uibar\rangle  \notag \\
&=&\sum_{m\geq 0}\langle a,u_i(m-1)b(2l-m-1)\uibar\rangle
\xi ^{m}.  \quad \label{Fxy_expan2}
\end{eqnarray}%
Since $l$ is the  lowest primary weight, we have $b(2l-m-1)\uibar\in V_{m}=\left(V_{\omega}\right)_m$ for $0\le m<l$  which  determines the first $l$ coefficients in the $\xi$ expansion \eqref{Fxy_expan2}.
This follows by writing  $b(2l-m-1)\uibar$ in a Virasoro basis  with  coefficients computed in a similar way as for the Casimir vectors in Lemma~\ref{lemma:lambda}.
On the other hand, from \eqref{FabG} we find using $y=-\xi x/(1-\xi)$ that 
\begin{eqnarray*}
y^{2l}F(a,b;x,y)
&=&
g\left(-\frac{\xi }{1-\xi}\right)\left(1-\xi\right)^{2l}\\
&=& 
A_{{0}}- \left( 2lA_{{0}}+A_{{1}} \right) \xi +O(\xi^2),
\end{eqnarray*}
for  $g(y)=G(a,b;1,y)=\sum_{r=0}^{4l}A_{r}y^{r}$. 
Hence the first $l$ coefficients of \eqref{Fxy_expan2} determine $A_{0},\ldots , A_{l-1}$.
Thus, using $b(2l-1)\uibar=(-1)^l \langle b,\uibar\rangle \vac$, we have
\[
A_{0}=\langle a,u_{i}(-1)b(2l-1)\uibar\rangle
=(-1)^l \langle a,\ui\rangle \langle b,\uibar\rangle =(-1)^l \langle
a,b\rangle.
\]
In general, $A_k=\langle a,b\rangle a_k(c)$ for $k=0,\ldots, l-1$ for a rational function $a_k(c)$.
The other $l+1$ coefficients of $g(y)$  (recalling $A_{r}=A_{4l-r}$) are determined by using associativity \eqref{Assoc} and expanding in $\zeta=(x-y)/y$ as follows
\begin{eqnarray}
\left(x-y\right)^{2l}F(a,b;x,y)
&=&
\sum_{m\in \Z}
\langle a,Y\left(\ui(m) \uibar,y\right)b\rangle
 (x-y)^{2l-m-1}
\notag
\\
&=&
\sum_{n\ge 0} B_n \zeta^{n},
\label{FB}
\end{eqnarray}
for $B_n=\langle a,o(\lambda^{(n)})b\rangle$ for $n\ge 0$ and recalling $o(\lambda^{(n)})=\lambda^{(n)}(n-1)$.
\begin{lemma}
\label{lem:Bn}
The leading coefficients of \eqref{FB} are $
B_0=(-1)^l p_l \langle a, b\rangle $ 
and $B_1=0$. 
For $k\ge 1$, the odd labelled coefficients $B_{2k+1}$ obey 
\begin{eqnarray*}
B_{2k+1}=\frac{1}{2}\sum_{r=2}^{2k}\binom{-r}{2k+1-r}(-1)^r B_r,
\end{eqnarray*}
i.e. $B_{2k+1}$ is determined by the lower even labelled coefficients $B_2,\ldots ,B_{2k}$.
The even labelled coefficients are given for $k\ge 0$ by 
\begin{eqnarray}
B_{2k}&=&A_{2l}\delta_{k,0}+\sum_{m=1}^{2l} \left[ \binom{m}{2k}+ \binom{-m}{2k}\right]A_{2l-m}. \label{B2k}
\end{eqnarray}
\end{lemma}
\noindent \textbf{Proof.}
From Lemma~\ref{lemma:lambda01} we have $\lambda^{(0)}= (-1)^l p_l\vac$
and $\lambda^{(1)}=0$ so that $B_0=(-1)^l p_l \langle a, b\rangle$ and $B_1=0$.
Comparing \eqref{FB} to \eqref{FabG} we find that 
\begin{eqnarray*}
\sum_{n\ge 0} B_n \zeta^{n}
&=&
g\left(\frac{1}{1+\zeta}\right)\left(1+\zeta\right)^{2l}=g\left(1+\zeta\right)\left(1+\zeta\right)^{-2l},
\end{eqnarray*}
since $G(a,b;x,y)$ is symmetric and homogeneous. Thus 
\[
\sum_{n\ge 0} B_n \zeta^{n}=\sum_{n\ge 0} B_n \left(\frac{-\zeta}{1+\zeta}\right)^{n}.
\]
This implies $B_n=\sum_{r=0}^{n}\binom{-r}{n-r}(-1)^r B_r$. Taking $n=2k+1$ leads to the stated result. 
\eqref{B2k} follows from the identity
\begin{eqnarray*}
\sum_{n\ge 0} B_n \zeta^{n}
&=& A_{2l}+\sum_{m=1}^{2l} A_{2l-m}\left[ (1+\zeta)^m+(1+\zeta)^{-m}\right].\quad \qed
\end{eqnarray*}
\
We next assume that $\lambda^{(n)}\in V_{\omega}$ for even $n\le 2l$ giving $B_{2k}=p_l\langle a, b\rangle b_{2k}(c)$ for $k=1,\ldots,l$ for some rational functions $b_{2k}(c)$  via Lemma~\ref{lem:ratlam}. 
Note that we are not (yet) assuming  $\lambda^{(2l+2)}\in V_{\omega} $.
$G(a,b;x,y)$ is uniquely  determined provided we can invert \eqref{B2k} to solve for $A_{l},\ldots , A_{2l}$. Define the $l\times l$ matrix  
\begin{eqnarray}
\label{eq:V}
M_{mk}= \binom{m}{2k}+ \binom{-m}{2k},
\end{eqnarray}
of coefficients for $A_{2l-m}$ of $B_{2k}$ in \eqref{B2k},
where $m,k=1,\ldots, l$.
\begin{lemma}  $M$ is invertible with
 $\det M=1$.
\end{lemma}
\noindent \textbf{Proof.}
Define  unit diagonal lower and upper triangular matrices $L$ and $U$ by
\begin{eqnarray*}
L_{ij}=\left\{
\begin{array}{cc}
	\binom{2i-j-1}{j-1}& \mbox{for } i\le j ,\\
	0 & \mbox{for } i>j,
\end{array}
\right.
\qquad
U_{jk}=
\left\{
\begin{array}{cc}
	\frac{k}{j}\binom{j+k-1}{2j-1}  &\mbox{for } j \le k,\\
	0 & \mbox{for }j>k.
\end{array}
\right.
\end{eqnarray*}
By induction in $k$, one can show that $M_{ik}=(L\, U)_{ik}$ and so $\det M=1$. \tqed
\
Next assume $\lambda^{(2l+2)}\in V_{\omega} $ giving another condition on $B_{2l+2}$ (already determined from \eqref{B2k}). Thus $p_l$ is a specific rational function of $c$. Hence we have
\begin{proposition}
\label{prop:Fab}
Let $V$ be an Exceptional VOA with lowest primary weight $l$. Then the genus zero correlation function $F(a,b;x,y)$ is uniquely determined and $p_l=p_l(c)$, a specific rational function of $c$. 
\end{proposition}
For $l=1,2$ we may use $F(a,b;x,y)$ to understand many properties of the corresponding VOA (as briefly reviewed below) \cite{T1,T2}. 
We already know from Proposition~\ref{prop:Zgen}(ii) that $p_l=\dim V_l-\dim \left(V_{\omega}\right)_l$ is a rational function of $c$. 
In principle, the specific rational expressions for $p_l$ may differ but, in practice, the same expression is observed to arise for all $l\le 9$. A more significant point is that the argument leading to  Proposition~\ref{prop:Fab} may be adopted to understanding some automorphism group properties of $V$.

\subsection{Exceptional VOAs of Class ${\mathcal S}^{2l+2}$}
Let $G=\Aut(V)$ denote the automorphism group of a VOA $V$ and let $V^G$ denote the sub--VOA fixed by $G$.
Since the Virasoro vector is $G$ invariant it follows that $V_{\omega}\subseteq V^G$.
$V$ is said to be of Class ${\mathcal S}^n$ if $V_k^G=\left(V_{\omega}\right)_k$ for all $k\le n$ \cite{Mat}. 
(The  related notion of conformal $t$--designs is described in \cite{Ho2}.)
In particular, the quadratic Casimir \eqref{eq:lamdef} is $G$--invariant so it follows that a VOA $V$ with lowest primary weight $l$ of class ${\mathcal S}^{2l+2}$ is an Exceptional VOA. It is not known if every Exceptional VOA is of class ${\mathcal S}^{2l+2}$.

The primary vector space $\Pi_l$ is a finite dimensional $G$--module. Assuming  $\Pi_l$ is a reducible $G$--module (e.g.  for $G$ linearly reductive  \cite{Sp}) we have
\begin{proposition}
\label{prop:irred}
Let $V$ be an Exceptional VOA of class $\mathcal{S}^{2l+2}$ with primaries ${\Pi_l}$ of  lowest weight $l$. If $\Pi_l$ is a reducible $G$--module then it is either an irreducible $G$--module or the direct sum of two isomorphic irreducible $G$--modules.
\end{proposition}
\begin{remark}
\label{rem:oddpl}
For odd  $p_l$  it follows that $\Pi_l$ must be an irreducible $G$--module.
\end{remark}

\noindent \textbf{Proof.}
Let $\rho$ be a $G$--irreducible component of $\Pi_l$ and let $\overline{\rho}$ denote the $\langle\ ,\ \rangle$  dual vector space.
$\overline{\rho}$ and $\rho$ are isomorphic as $G$--modules.
Define 
\begin{eqnarray*}
R=\left\lbrace
\begin{array}{ll}
	\rho & \mbox{if } \rho=\overline{\rho},\\
	\rho\oplus \overline{\rho}& \mbox{if } \rho\neq\overline{\rho}.
\end{array}
 \right.
\label{eq:Rdef}
\end{eqnarray*}
Clearly $R\subseteq \Pi_l$ is a self-dual vector space. 
We next repeat the Casimir construction and  analysis that lead up to  Proposition~\ref{prop:Fab}.
Choose an $R$--basis $\{\vi:i=1,\ldots ,\dim R\}$  and dual basis $\{\vibar\}$ and define Casimir vectors 
\begin{eqnarray}
\lam{(n)}_R=\vi(2l-n-1)\vibar\in V_{n},\quad n\ge 0,
\label{eq:lamdefR}
\end{eqnarray}
where now we sum $i$ from 1 to $\dim R\le p_l$. 
But $\lam{(n)}_R$ is $G$--invariant and since $V$ is of class ${\mathcal S}^{2l+2}$, it follows that 
$\lam{(n)}_R\in V_{\omega}$ for all $n\le 2l+2$. We define a genus zero correlation function
constructed from the vector space $R$
\begin{eqnarray}
F_R(a,b;x,y)=\langle a,Y(\vi,x)Y(\vibar,y)b\rangle,\label{eq:FabR}
\end{eqnarray}
for all $a,b\in R$. We then repeat the earlier arguments to conclude that  Proposition~\ref{prop:Fab} also holds for $F_R(a,b;x,y)$ where, in particular, $\dim R=p_l(c)$, for the \textbf{same} rational function. Thus  $\dim R=l$ and the result follows.
\tqed
\
%
 %
 %

\section{Exceptional VOAs of Lowest Primary Weight $l\le 9$}
\label{sect: VOAex}
We now consider Exceptional VOAs of lowest primary weight $l\le 9$. We  denote  
by $E_n=E_n(q)$ the Eisenstein series of weight $n$ appearing in the MLDE \eqref{eq:MLDE}. 
For $l\le 4$ we describe all the rational values for $c,h$ whereas for $5\le l \le 9$ we give all rational values for $c,h$ for which $p_l=\dim \Pi_l\le 500000$,
found by computer algebra techniques. We also consider  conjectured extremal self-dual VOAs with $c=24(l-1)$ \cite{Ho1, Wi}.
Any MLDE solution for rational $h$ for which there is no irreducible character is marked with an asterisk. 
We obtain many examples of known Exceptional VOAs such Deligne's Exceptional Series of Lie algebras, the Moonshine and Baby Monster modules. 
There are also a number of candidate solutions for which no construction yet exists indicated  by question marks. 
\medskip

\noindent $\boldsymbol{[l=1].}$
This  is discussed in much greater detail in \cite{T1,T2}.  Propositions~\ref{prop:MLDE}--\ref{prop:Zgen}
imply that $Z_{N}(q)$ satisfies the following 2nd order MLDE \cite{T2}
\begin{eqnarray*}
D^2 Z-\frac{5}{4}c \left( c+4 \right) E_{4}Z =0.
\end{eqnarray*}
This MLDE has also appeared  in 
\cite{MatMS, KZ, Mas2, KKS, Kaw}. 
The indicial roots $x_1=-c/24, x_2= (c+4)/24$ are exchanged under the MLDE symmetry $c\leftrightarrow -c-24$. Solving iteratively for the partition function
\begin{eqnarray*}
Z_{V}(q)=q^{-c/24}\left(1+p_1q+(1+p_1+p_2)q^2+(1+2p_1+p_2+p_3)q^3+\ldots\right),
\end{eqnarray*}
where $p_n=\dim \Pi_n$, for weight $n$ primary vector space $\Pi_n$, we have
\begin{eqnarray*}
&&p_1=
\frac {c ( 5 c +22 ) }{10-c},
\quad p_{2}=
{\frac {5 ( 5 c +22 )  ( c-1 )  ( c+2) ^{2}}{2 ( c-10 )  ( c-22 ) }},\\
&&p_{3}=-{\frac {5c ( 5c + 22 )  ( c-1 )  ( c+5 )  ( 5 {c}^{2}+268 ) }{ 6( c-10 ) 
 ( c-22 )  ( c-34 ) }},\ldots
\end{eqnarray*}
For $c=10\mod 12$, the indicial roots differ by an integer leading to denominator zeros for all $p_n$.

By Proposition~\ref{prop:Zgen}, $V$ is generated by $V_1$ which defines a Lie algebra $\frak{g}$. 
$F(a,b;x,y)$ from Proposition~\ref{prop:Fab} determines the Killing form which can be used to show that $\frak{g}$ is simple with dual Coxeter number \cite{T1,MT}
\[
h^{\vee }=6k\frac{2+c}{10-c},
\] 
for some real level $k$. Thus
$V=V_{\mathfrak{g}}(k)$, a level $k$ Kac-Moody VOA.

The  indicial root $x_2$ of the MLDE gives the lowest weight $h=(c+2)/12$ of any independent irreducible $V$--module(s) $N$. 
Therefore $V_{\mathfrak{g}}(k)$ has at most two independent irreducible characters so that the level $k$ must be positive integral \cite{Kac2}. 
Comparing $p_1$ and $h^{\vee }$ to Cartan's list of simple Lie algebras shows that in fact $k=1$ with  $c=1,2,\frac{14}{5},4,\frac{26}{5},6,7,8$ with
$
\mathfrak{g}=A_1,A_2,G_2,D_4,F_4,E_6,E_7,E_8,
$
respectively, 
known as the Deligne Exceptional Series \cite{D,DdeM,MarMS,T2}. 
In summary, we have
\begin{center}
\renewcommand{\arraystretch}{1.5}
\begin{tabular}{|c|c |c|c|c|c|}
\hline
$c>0$ & $p_{1}$ & $p_2$ & $p_3$ &VOA & $h \in \Q$\\
\hline
$1$ & 3 & 0 & 0 & $V_{A_1}(1)$ & $\frac{1}{4}$\\
\hline
$2$ & 8 &  8 & 21 & $V_{A_2}(1)$& $\frac{1}{3}$\\
\hline
$\frac{14}{5}$ & 14 &  27 & 84 & $V_{G_2}(1)$& $\frac{2}{5}$\\
\hline
$4$ & 28 & 105 & 406 &  $V_{D_4}(1)$& $\frac{1}{2}$\\
\hline
$\frac{26}{5}$ & 52 & 324 & 1547 &  $V_{F_4}(1)$& $\frac{3}{5}$\\
\hline
$6$ & 78 & 650 & 3575 &  $V_{E_6}(1)$ & $\frac{2}{3}$\\
\hline
$7$ & 133 & 1539 & 10108 &  $V_{E_7}(1)$& $\frac{3}{4}$\\
\hline
$8$ & 248 & 3875 & 30380 & $V_{E_8}(1)$& $\frac{5}{6}*$\\
\hline
\end{tabular}
\end{center}
The table also shows $h$ for a possible  irreducible $V$-module(s). For $c=2$ and $4$ there are 2 independent irreducible modules but which share the same character (due to $\mathfrak g$ outer automorphisms ). 
$V_{E_8}(1)$ is self-dual so that the MLDE solution with $h=\frac{5}{6}$ is not an irreducible character.
\medskip

\noindent $\boldsymbol{[l=2].}$
This case is also discussed in detail in \cite{Mat,T1,T2}.  Propositions~\ref{prop:MLDE}--\ref{prop:Zgen}
imply that $Z_{N}(q)$ satisfies the following 3rd order MLDE \cite{T2}
\begin{eqnarray*}
&&D^3 Z-{\frac {5}{124}} \left( 704+240 c+21 {c}^{2}\right) E_{4}\, D Z
-{\frac {35}{248}} c \left( 144+66 c+5{c}^{2} \right)E_{{6}}\, Z=0,
\end{eqnarray*} 
with indicial equation \eqref{eq:indicial}
\[
(x-x_1)\left(x^2-\left(\half+x_1\right)x
+\frac{20x_{1}^{2}-11x_{1}+1}{62}\right)=0,
\]
for $x_1=-c/24$. 
Solving iteratively for the partition function ($x=x_1$)
\begin{eqnarray*}
Z_V(q)=q^{-c/24}(1+(1+p_2){q}^{2}+(1+p_2+p_3){q}^{3}+\ldots),
\end{eqnarray*}
where $p_n=\dim \Pi_n$, for weight $n$ primary vector space $\Pi_n$, we find
\begin{eqnarray*}
p_2=
\frac { ( 7 c +68)  ( 2 c-1 )  ( 5 c +22) }
{2(c^2-55 c+748)},\quad 
p_3=
\frac {31c ( 7 c+ 68)  ( 2 c-1 )  ( 5 c+44 )  ( 5 c+22 ) }
{ 6( c^2-55 c+748 )  ( c^2-86 c+1864 ) }.
\end{eqnarray*}
From Proposition~\ref{prop:Zgen}, the Griess algebra generates $V$ and from Proposition~\ref{prop:Fab}  the Griess algebra is simple \cite{T1}. 
This leads to the following possible Exceptional VOAs with $c,h\in \Q$
\begin{center}
\renewcommand{\arraystretch}{1.7}
\begin{tabular}{|c|c|c|c|c|}
\hline
$c$ & $p_{2}$ & $p_{3}$& VOA& $h\in\Q$\\
\hline
$-\frac{44}{5}$ & 1 & $ 0$ & $L\left(c_{3,10},0\right)\oplus L\left(c_{3,10},2\right)$ & $0,-\frac{1}{5},-\frac{2}{5}$\\
\hline
$8$ & 155 & 868 &  $V_{\sqrt{2}E_8}^{+}$ & $0,\frac{1}{2},1$\\
\hline
$16$ & 2295 & $ 63240$ & $V_{\BW_{16}}^{+}$  & $0,1,\frac{3}{2}$\\
\hline
$\frac{47}{2}$ & 96255 & $ 9550635$ & $\VB ^{\natural}_{\Z}$ & $0,\frac{3}{2},\frac{31}{16}$\\
\hline
$24$ & 196883 & $ 21296876$ & $V^{\natural}$ & $0$\\
\hline
$32$ &  $3 . 7^{2} . 13 .  73 $ & $    2^{4} .  3 .    7^{2}  .  13 .  31 .73 $ & ?? $V_{L}^{+}\oplus (V_L)_T^{+}$; $L$ extremal & $0$\\
\hline
$\frac{164}{5}$ & $3^2.17.19.31$ &
 $ 2.5.13.17.19.31.41 $ 
& ??& $0,\frac{11}{5},\frac{12}{5}$\\
\hline
$\frac{236}{7}$ & $5.19.23.29$
 &  $2.19.23.29.31.59$ 
& ??& $0,\frac{16}{7},\frac{17}{7}$\\
\hline
$40$ & $3^2.29.79$ &
 $ 2^2.5 .29.31.61.79$ &?? $V_{L}^{+}\oplus (V_L)_T^{+}$; $L$ extremal& $0$\\
\hline
\end{tabular}
\end{center}
The list includes the famous Moonshine Module $V^{\natural}$ \cite{FLM}, the Baby Monster VOA $\VB ^{\natural}_{\Z}$ \cite{Ho1},   $V_{L}^{+}$  for $L=\sqrt{2}E_8$ \cite{G}  and the rank 16 Barnes-Wall lattice $L=\BW_{16}$ \cite{Sh}, and a minimal model simple current extension $AD$--type as in Proposition~\ref{prop:ADVOA}. 
The value(s) of $h=x_i+c/24$  for the lowest weight(s) agree with those for the irreducible $V$-modules as do the corresponding MLDE solutions for the characters in each case.
There are also four other possible candidates. For $c=32$ and 40 one can construct a self-dual VOA from an extremal even self-dual lattice $L$ (with no vectors squared length 2).  
However, such lattices are not unique and it is not known which, if any, gives rise to a VOA satisfying the exceptional conditions. There are no known candidate constructions for $c=\frac{164}{5}$ and $\frac{236}{7}$.  

Note that $p_2=\dim\Pi_2$ is odd in every case and Proposition~\ref{prop:irred} implies that if $\Pi_2$ is $\Aut V$--reducible then it is irreducible. This is indeed the case in the first five known cases for $c\le 24$ \cite{Atlas}. 
$\Pi_3$ is also an $\Aut V$-module whose dimension $p_3$ is  given. 
The MLDE solutions (with positive coprime integer coefficients) for $c=164/5$ with $h=11/5, 12/5$ 
and for $c=236/7$ with $h=16/7, 17/7$  have respective leading $q$--expansions
\begin{eqnarray*}
Z_{11/5}(q)&=&q^{5/6}\left(
2^3.31.41
\,+\,
5.11.31.41.53
\,q+O(q^2)
\right),\\
Z_{12/5}(q)&=&q^{31/30}\left(
2^2.11^2.31.41
\,+\,
2^5.11^2.31^2.41
\,q+O(q^2)\right),\\
Z_{16/7}(q)&=&q^{37/42}\left(
17.23.31 
 \,+\,   
2^5.7.17.31.37 
\,q+O(q^2)
\right),\\
Z_{17/7}(q)&=&q^{43/42}\left(
2^4.29.31.59
\,+\,
2.3.17.29.31.43.59
\,q+O(q^2)\right).
\end{eqnarray*}
These coefficients constrain the possible  structure of  $\Aut V$ further.

\medskip

\noindent $\boldsymbol{[l=3].}$
$Z_{N}(q)$ satisfies the 4th order MLDE  
\begin{eqnarray*}
&& \left( 578 c -7 \right) \,D^4  Z
 -\frac{5}{2}  \left( 168 {c }^{3}+2979 {c }^{2}+15884 c -4936 \right) E_{4}\,D^2  Z
\\
&&-{\frac {35}{2}}\left( 25 {c }^{4}+661 {c }^{3}+4368 {c }^{2}+10852 c +1144
 \right) E_{{6}}\,D Z
\\
&&
-{\frac {75}{16}} c  \left( 14 {c }^{4}+
425 {c }^{3}+3672 {c }^{2}+5568 c +9216 \right) E_{4}^{2}\, Z=0.
\end{eqnarray*}
Solving iteratively for the partition function
we find \cite{T2}
\begin{eqnarray*}
Z_V&=&q^{-c/24}(1+{q}^{2}+(1+p_3){q}^{3}+(2+p_3+p_4){q}^{3}+\ldots),\\
p_3&=&-{\frac {(5c+22)(3c+46)(2c-1)(5c+3)(7c+68) }{5 {c}^{4}-
703 {c}^{3}+32992 {c}^{2}-517172 c+3984}},
\end{eqnarray*}
and $p_4=\frac{r(c)}{s(c)}$ for 
\begin{eqnarray*}
r(c)&=& 
-\frac{1}{2} 
\left( 2 c-1 \right)  \left( 3 c+46 \right)  \left( 5
 c-4 \right)  \left( 7 c+68 \right)  \left( 5 c+3 \right)  \left( 7 c+114 \right) \\
&&.  \left( 55 {c}^{3}-5148 {c}^{2}-11980 c-36528\right) ,
\\
s(c)&=& 
\left( 5 {c}^{4}-703 {c}^{3}+32992 {c}^{2}-517172 c+3984 \right)
\\
&& . \left( 5 {c}^{4}-964 {c}^{3}+62392 {c}^{2}-1355672 c
+13344 \right).
\end{eqnarray*}
The $c,h\in \Q$ solutions  for positive integer $p_3$ with possible VOAs are
\begin{center}
\renewcommand{\arraystretch}{1.7}
\begin{tabular}{|c|c|c |c|c|}
\hline
$c$ & $p_{3}$ & $p_{4}$ & VOA & $h\in \Q$\\
\hline
$-{\frac {114}{7}}$ & 1 & 0 & $L\left(c_{3,14},0\right)\oplus L\left(c_{3,14},3\right)$ & $0,-\frac{3}{7}, -\frac{4}{7},-\frac{5}{7}$\\
\hline
${\frac {4}{5}}$ & 1 &  0 & $L\left(c_{5, 6},0\right)\oplus L\left(c_{5, 6},3\right)$ & $0,\frac{1}{15},\frac{2}{5},\frac{2}{3}$\\
\hline
$48$ & $3^2.19^2.101.131$ 
& $5^6.19^2.71.101$ 
&?? H\"ohn  Extremal VOA & $0$ \\
\hline 
\end{tabular}
\end{center}
The H\"ohn  Extremal VOA is a conjectural self-dual VOA \cite{Ho1}. If $\Pi_3$ is a reducible $\Aut(V)$--module then it must be irreducible excluding Witten's suggestion that $\Aut(V)=\mathbb{M}$, the Monster group  \cite{Wi}.   

\medskip
\noindent $\boldsymbol{[l=4].}$
Proposition~\ref{prop:MLDE} implies $Z_{N}(q)$  satisfies the 5rd order MLDE 
\begin{eqnarray*}
&& \left( 317 c+3 \right)  \, D^{5}Z
-\frac{5}{7} \left( 297 {c}^{3}+6746 {c}^{2}+53133 c+4644 \right)   E_{{4}}\, D^{3}Z
\\
&&
-{\frac {25}{8}}  \left( 77 {c}^{4}+3057 {c}^{3}+31506 {c}^{2}+129736 c-
24096 \right)  E_{{6}}\, D^{2}Z
\\
&&
-{\frac {25}{112}}  \left( 231 {c}^{5}+12117 {c}^{4}+194916 {c}^{3}+843728 {c}^{2}+1652288 c-718080 \right) {E_{{4}}}^{2} \,D Z
\\
&&
-{\frac {25}{32}} c\left( c+24 \right)  \left( 15 {c}^{4}+527 {c}^{3}+5786 {c}^{2}+528 c+25344 \right) E_{{4}}E_{{6}} \,Z=0.
\end{eqnarray*}
Solving iteratively for the partition function
we find 
\begin{eqnarray*}
Z_V&=&q^{-c/24}(1+{q}^{2}+{q}^{3}+(2+p_4)q^4+(3+p_4+p_5)q^5+\ldots),\\
p_{4} &=&
\frac { 5 ( 3 c+46  )   ( 2 c-1  )   ( 11 c+232  )   ( 7 c + 68 )   ( 5 c+3  )   ( c+10  ) }
{2  ( 5 c^4-1006 c^3+67966 c^2-1542764 c-12576  )   ( c-67  ) },
\end{eqnarray*}
and $p_5=\frac{r(c)}{s(c)}$ where 
\begin{eqnarray*}
r(c)&=&
3( c-1 )  ( 5c +22)  ( 3c+46 )  ( 2c-1 )  ( 11c+232 )  ( 7c+68 )  ( 5c+3 )  ( c+24 )
\\
&&  .( 59c^3-13554c^2+788182c-398640 ) 
,\\
s(c)&=&
2 ( c-67 )  ( 5c^4-1006c^3+67966c^2-1542764c-12576 )
\\
&& . ( 5c^{5}-1713c^4+221398c^3-
12792006c^2+278704260c+2426976 ) 
.
\end{eqnarray*}
The $c,h\in \Q$ solutions  for $p_4\le 500000$ and $c=48$ with possible VOAs are
\begin{center}
\renewcommand{\arraystretch}{1.7}
\begin{tabular}{|c|c|c |c|c|}
\hline
$c$ & $p_{4}$ &$p_{5}$ & VOA& $h\in \Q$\\
\hline
$1$ & $2$ & 
$0
$ & 
$V_{L}^{+}$ for $L=2\sqrt{2}\Z$& $0,\frac{1}{16},\frac{1}{4},\frac{9}{16},1$  \\
\hline
$72$ & $2^3. 11^4. 13^2. 131$ & 
$2. 11^4. 13^2. 103. 131. 191$ & 
?? H\"ohn Extremal  VOA &0 \\
\hline
\end{tabular}
\end{center}
The H\"ohn  Extremal VOA is a conjectured self-dual VOA \cite{Ho1,Wi}.
 If $\Pi_4$ is a reducible $\Aut(V)$--module,  then by Proposition~\ref{prop:irred}, either $p_4$ or $\frac{1}{2}p_4$ is the dimension of an irreducible $\Aut(V)$--module.
 
\medskip
\noindent $\boldsymbol{[l=5].}$
$Z_V$ satisfies a 6th order MLDE 
 with $p_5=\frac{r(c)}{s(c)}$ for 
\begin{eqnarray*}
r(c)&=&-( 13 c+350 )  ( 7 c+25 )  ( 5 c+126 )  ( 11 c+232 )\\
&& . ( 2 c-1 )  ( 3 c+46 )  ( 68+7 c )  ( 5 c+3 )  ( 10 c-7 ),\\
s(c)&=&1750 c^8-760575 c^7+132180881 c^6-11429170478 c^5\\
&&+484484459322 c^4-7407871790404 c^3-37323519053016 c^2\\
&&+25483483057200 c-363772080000.
\end{eqnarray*}
The $c,h\in \Q$ solutions  for $p_5\le 500000$ with possible VOAs are
\begin{center}
\renewcommand{\arraystretch}{1.7}
\begin{tabular}{|c|c |c|c|}
\hline
$c$ & $p_{5}$ & VOA & $h\in \Q$\\
\hline
$-\frac{350}{11}$ & 1 & $L\left(c_{3,22},0\right)\oplus L\left(c_{3,22},5\right)$
 & $0,-{\frac {8}{11}},-{\frac {10}{11}},-{\frac {13}{11}},-{\frac {14}{11}},-{\frac {15}{11}}
$\\
\hline
$\frac{6}{7}$ & 1 & $L\left(c_{6,7},0\right)\oplus L\left(c_{6,7},5\right)$ &
$0,{\frac {1}{21}},{\frac {1}{7}},{\frac {10}{21}},{\frac {5}{7}},{\frac {4}{3}}$\\
\hline
\end{tabular}
\end{center}
Witten's conjectured Extremal VOA for $c=4.24=96$ does not appear \cite{Wi}.

\medskip
\noindent $\boldsymbol{[l=6].}$
$Z_V$ satisfies a 7th order MLDE 
 with $p_6=\frac{r(c)}{s(c)}$ for 
\begin{eqnarray*}
r(c)&=&\frac {7}{2} ( 13 c +350 )  ( 5 c +164 )  ( 7 c +25 )  ( 11 c +232 )  ( 3 c +46 )
\\
&& . ( 4 c +21 )  ( 5 c +3 )  ( 10 c -7 )  ( 5 c ^2+316 c +3600 ) 
,\\
s(c)&=&1750 c ^9-1119950 c ^8+297661895 c ^7-41808629963 c ^6
\\
&&+3225664221176 c ^5-123384054679580 c ^4+1266443996541232 c ^3
\\
&&+29763510364647840 c ^2+96385155929078400 c +7743915615744000.
\end{eqnarray*}
The $c,h\in \Q$ solutions  for $p_6\le 500000$ with possible VOAs are
\begin{center}
\renewcommand{\arraystretch}{1.7}
\begin{tabular}{|c|c |c|c|}
\hline
$c$ & $p_{6}$ & VOA & $h\in \Q$\\
\hline
$-\frac{516}{13}$ & 1 & $L\left(c_{3,26},0\right)\oplus L\left(c_{3,26},6\right)$ &
$0,-{\frac {10}{13}},-{\frac {15}{13}},-{\frac {17}{13}}$\\
 &  &  & $-{\frac {20}{13}},-{\frac {21}{13}},-{\frac {22}{13}}$ \\
\hline
$-47$ & 1 & ${\mathcal W}(6)$ & $0,-\frac{5}{4},-\frac{3}{2},-\frac{5}{3}$ \\
 &  &  & $-{\frac {15}{8}},-{\frac {23}{12}},-2$
\\
\hline
$120$ & $2 .7^2 .11. 29 .43 .67.97 . 191 $& ?? Witten  Extremal  VOA & 0\\
\hline
\end{tabular}
\end{center}
$c=-47$ is first  example of a $\mathcal W(3k)$--algebra of Proposition~\ref{prop:W3kVOA}. The irreducible lowest weight $h$ values  and  character solutions agree with  \cite{F}.
Witten's conjecture Extremal VOA for $c=5.24=120$ appears \cite{Wi} where either $p_6$ or $\frac{1}{2}p_6$ is the dimension of an irreducible $\Aut(V)$--module.

\medskip

\noindent $\boldsymbol{[l=7].}$
$Z_V$ satisfies an 8th order MLDE 
 where $p_7=\frac{r(c)}{s(c)}$ for 
\begin{eqnarray*}
r(c)&=&-5 
\left( 13 c+350 \right)  \left( 5 c+164 \right) 
 \left( 7 c+25 \right)  \left( 11 c+232 \right)  \left( 3 c+46
 \right)  \left( 17 c+658 \right)
\\
&&  .\left( 4 c+21 \right)  \left( 5 
c+3 \right)  \left( 10 c-7 \right)  \left( 35 {c}^{3}+3750 {c}^{2}+
76744 c-32640 \right) 
,\\
s(c)&=&61250 {c}^{11}-54725125 {c}^{10}+
20922774275 {c}^{9}-4421902106730 {c}^{8}
\\
&&
+553932117001488 {c}^{7}-
40395124111104312 {c}^{6}+1491080056338817984 {c}^{5}\\
&&
-12528046696953576896 {c}^{4}-483238055074755678656 {c}^{3}\\
&&
-1702959754355175160320 {c}^{2}+249488376255167616000 c
\\
&&+362620505915136000000.
\end{eqnarray*}
There are no rational $c$ solutions for $p_7\le 500000$.

\medskip

\noindent $\boldsymbol{[l=8].}$
$Z_{V}(q)$ satisfies a 9th order MLDE with one $c,h\in \Q$ solution for $p_8\le 500000$
\begin{center}
\renewcommand{\arraystretch}{1.7}
\begin{tabular}{|c|c |c|c|}
\hline
$c$ & $p_{8}$ & VOA & $h\in \Q$\\
\hline
$-\frac{944}{17}$ & 1 & $L\left(c_{3,34},0\right)\oplus L\left(c_{3,34},8\right)$ &
$0,-{\frac {14}{17}},-{\frac {25}{17}},-{\frac {28}{17}},-{\frac {33}{17}}$ \\
 &  &  & $-{\frac {35}{17}},-{\frac {38}{17}},-{\frac {39}{17}},-{\frac {40}{17}}$ \\
\hline
\end{tabular}
\end{center}
\medskip

\noindent $\boldsymbol{[l=9].}$
$Z_{V}(q)$ satisfies a 10th order MLDE with $c,h\in \Q$ solutions  for $p_9\le 500000$
\begin{center}
\renewcommand{\arraystretch}{1.7}
\begin{tabular}{|c|c |c|c|}
\hline
$c$ & $p_{9}$ & VOA & $h\in \Q$\\
\hline
$-{\frac {1206}{19}}$ & 1 & $L\left(c_{3,38},0\right)\oplus L\left(c_{3,38},9\right)$ &
$0,-{\frac {16}{19}},-{\frac {29}{19}},-{\frac {36}{19}},-{\frac {44}{19}}$ \\
 &  &  & $-{\frac {46}{19}},
-{\frac {49}{19}},-{\frac {50}{19}},-{\frac {51}{19}}
$ \\
\hline
$-{\frac {208}{35}}$ & 1 & $L\left(c_{5,14},0\right)\oplus L\left(c_{5,14},9\right)$ &
$0,-\frac{2}{7},-{\frac {9}{35}},-{\frac {4}{35}}
$ \\
 &  &  & $\frac{1}{7},{\frac {11}{35}},{\frac {9}{7}},\frac{8}{5}
$ \\
\hline
$-{\frac {14}{11}}$ & 1 & $L\left(c_{6,11},0\right)\oplus L\left(c_{6,11},9\right)$ &
$0,-\frac{1}{11},-{\frac {2}{33}},\frac{1}{11},{\frac {7}{33}}
$ \\
 &  &  & 
$
{\frac {6}{11}},{\frac {25}{33}},{\frac {14}{11}},{\frac {52}{33}},\frac{8}{3}
$ \\
\hline
$-71$ & 1 & ${\mathcal W}(9)$ & $0,-2,-\frac{9}{4},-{\frac {39}{16}},-\frac{8}{3}
$ \\
 &  &  & 
$
-\frac{11}{4},-{\frac {35}{12}},-{\frac {47}{16}},-3,-{\frac {9}{8}}*
$ \\
\hline
\end{tabular}
\end{center}
For $c=-71$, the MLDE solutions  agree with all  the irreducible characters  for   ${\mathcal W}(9)$ in  \cite{F} except for $h=-{\frac {9}{8}}$.

\section{Exceptional VOSAs}
\subsection{VOSA Quadratic Casimirs and Zhu Theory}
We now give an analysis for Vertex Operator Superalgebras (VOSAs). Many of the results are similar but there are also significant differences e.g. here the MLDEs  involve twisted Eisenstein series. 
Let $V$ be a simple VOSA of strong CFT-type with unique invertible bilinear form $\langle\ ,\ \rangle$. 
Let $\Pi_l$ denote the space of Virasoro primary vectors of \emph{lowest half integer weight} $l\in \N+\half$ i.e. $\Pi_l$ is of odd parity and $V_k=(V_{\omega})_k$ for all $k\le l-\half$.
We  construct quadratic Casimir vectors $\lam{(n)}$ as in Section~\ref{sect:Cas} (from the odd parity space $\Pi_l$) which enjoy the same properties as VOA Casimir vectors.

Define the genus one partition function of a VOSA $V$ by
\begin{equation}
Z_V( q)=\Trq{V}{\sigma\,}=q^{-c/24}\sum_{n\geq 0}(-1)^{2n}\dim V_n \,q^n,\label{eq:Zsig(q)}
\end{equation}
for fermion number operator $\sigma$ where $\sigma u=(-1)^{p(u)}u$ for $u$ of parity $p(u)$ and with a corresponding definition for a simple ordinary $V$-module $N$. 
We also define the genus one 1-point correlation function  
\begin{equation}
Z_{N}(u,q)=\Trq{N}{\sigma\, o(u)}.
\label{eq:Zsig1pt}
\end{equation}
 In \cite{MTZ} a Zhu reduction formula for the 2-point correlation function $\Zq{N}{Y[u,z]v}$ for $u,v\in V$ is found expressed in terms of \textit{twisted elliptic Weierstrass functions} parameterized by $\theta,\phi\in\{\pm 1\}$. 
Let $\phi=e^{\tpi \kappa}$ for $\kappa\in \{0,\half\}$. 
Then \eqref{eq:Pm} and \eqref{eq:En} are generalized to \cite{MTZ}
\begin{equation}
P_{m}\begin{bmatrix}
\theta \\ 
\phi
\end{bmatrix} (z )
=\frac{1}{z^{m}}+(-1)^{m}\sum\limits_{n\geq m}\binom{n-1}{m-1}
E_{n}
\begin{bmatrix}
\theta \\ 
\phi
\end{bmatrix}
(q)z^{n-m},  \label{PkEntwist}
\end{equation}
for twisted Eisenstein series $E_{n}
\left[\begin{smallmatrix}
\theta \\ 
\phi
\end{smallmatrix}\right] (q )=0$  for $n$   odd,  
and for $n$  even 
\begin{eqnarray}
E_{n}\begin{bmatrix}
\theta \\ 
\phi
\end{bmatrix} (q ) &=&-\frac{B_{n}(\kappa )}{n!}+\frac{2}{(n-1)!}%
\sum\limits_{k\in\N+\kappa}\frac{k^{n-1}\theta q^{k}}{1-\theta q^{k}},  \label{Ektwist}
\end{eqnarray}%
and where $B_{n}(\kappa )$ is the Bernoulli polynomial defined by%
\begin{equation}
\frac{e^{z \kappa }}{e^z-1}=\frac{1}{z}+\sum\limits_{n\geq 1}\frac{B_{n}(\kappa )}{n!}z^{n-1}.  \label{Bernoulli}
\end{equation}%
\eqref{PkEntwist} and \eqref{Ektwist} agree with  \eqref{eq:Pm} and \eqref{eq:En} respectively for $(\theta,\phi)=(1,1)$.
$P_{m}
\left[\begin{smallmatrix}
\theta \\ 
\phi
\end{smallmatrix}\right] (z )$ converges absolutely and uniformly on compact subsets of the domain
$|q| < |e^{z} | < 1$ 
and $E_{n}
\left[\begin{smallmatrix}
\theta \\ 
\phi
\end{smallmatrix}\right] (q )$ is a holomorphic function of $q^{\half}$ for $|q|<1$.
 For $(\theta ,\phi )\neq (1,1)$, $E_{n}\left[\begin{smallmatrix}
\theta \\ 
\phi
\end{smallmatrix}\right] $ is modular of weight $n$ in the sense that 
\begin{equation}
E_{n}
\begin{bmatrix}
\theta ^{a}\phi ^{b} \\ 
\theta ^{c}\phi ^{d}
\end{bmatrix}
\left(\frac{a\tau
+b}{\gamma\tau +\delta}\right)=(\gamma\tau +\delta)^{n}E_{n}\begin{bmatrix}
\theta \\ 
\phi
\end{bmatrix}(\tau ).
\label{Ekmodtwist}
\end{equation}
Then Zhu reduction of Proposition~\ref{prop:Zhu} generalizes to \cite{MTZ}
\begin{proposition}
\label{prop:Zhutwist}
Let $N$ be a simple ordinary $V$--module for a VOSA $V$. For $u$ of parity $p(u)$ and for all $v$ we have
\begin{eqnarray*}
\Zq{N}{ Y[u,z]v}
&=&
\Trq{N}{\sigma\, o(u)o(v)}\delta_{p(u) 1}\notag \notag\\
&&+\sum_{m\ge 0} P_{m+1}
\begin{bmatrix}
1 \\ 
p(u)
\end{bmatrix}
(z) \Zq{N}{u[m]v}.
\end{eqnarray*}
\end{proposition}
For even parity $u$  this agrees with Proposition~\ref{prop:Zhu}.  In particular, 
Corollary~\eqref{thm:ZhuLk} concerning Virasoro vacuum descendents holds.
Much as before, Proposition~\ref{prop:Zhutwist} implies that the Casimir vectors $\lam{[n]}\in V_{[n]}$ obey
\begin{eqnarray}
\sum_{n\ge 0}\Zq{N}{\lam{[n]}}z^{n-2l}&=& 
\sum_{m=0}^{2l-1} P_{m+1}
\begin{bmatrix}
1 \\ 
-1
\end{bmatrix}
(z) \Zq{N}{\lam{[2l-m-1]}}
.\quad
\label{eq:lamrec}
\end{eqnarray}
Equating the  $z$ coefficients implies the following variant of Proposition~\ref{prop:lamrec}
\begin{proposition}
\label{prop:lamrectwist}
$\Zq{N}{\lam{[2l+1]}}$ satisfies the recursive identity
\begin{eqnarray}
\Zq{N}{\lam{[2l+1]}}&=&
-2\sum_{r=0}^{l-\half} 
(l-r) E_{2(l-r)+1}
\begin{bmatrix}
1 \\ 
-1
\end{bmatrix}
(q)
\Zq{N}{\lam{[2r]}}.\quad 
\label{Zrec2twist}
\end{eqnarray}
\end{proposition}

\subsection{Exceptional VOSAs} 
\label{sect:EVOSA}
Let $V$ be a simple VOSA of strong CFT--type  with primary vectors  of lowest weight $l\in \N+\half$ for which $\lam{(2l+1)}\in V_{\omega}$. 
We further assume that $ \left(V_{\omega}\right)_{2l+1}$ contains no Virasoro singular vectors. 
We call $V$ an \emph{Exceptional VOSA of Odd Parity Lowest Primary Weight $l$}. 
Proposition~\ref{prop:lam2k} implies
\begin{proposition}
\label{prop:MLDEtwist}
Let $V$ be an Exceptional VOSA of  lowest weight $l\in \N+\half$  and central charge $c$. Then  $Z_N(q)$  for a  simple ordinary $V$--module $N$  satisfies a Twisted Modular Linear Differential Equation (TMLDE) 
\begin{equation}
\sum_{m=0}^{l+\half} g_{l+\half-m}\begin{bmatrix}
1 \\ 
-1
\end{bmatrix}(q,c)\, D^{m}Z(q)=0,
\label{eq:TMLDE}
\end{equation}
where  
 $g_{k}\left[\begin{smallmatrix}
1 \\ 
-1
\end{smallmatrix}\right](q,c)$ is a twisted modular form of weight $2k$ whose coefficients over the ring of twisted Eisenstein series \eqref{Ektwist} are rational functions of $c$.
\end{proposition}

The TMLDE  \eqref{eq:TMLDE} is of order $l+\half$ with a regular singular point at $q=0$ provided $g_0\left[\begin{smallmatrix}
1 \\ 
-1
\end{smallmatrix}\right](q,c)=g_0(c)\neq 0$ so that Frobenius-Fuchs theory implies that its solutions are  holomorphic  in $q^{\half}$ for $0<|q|<1$. 
Furthermore, from \eqref{Ekmodtwist},  $\widehat{Z}_{N}=Z_{N}\left(\frac{\alpha\tau +\beta}{\gamma\tau +\delta}\right)$ is a  solution of the TMLDE
\begin{equation}
\sum_{m=0}^{l+\half} g_{l+\half-m}\begin{bmatrix}
(-1)^b \\ 
(-1)^c
\end{bmatrix}(q,c)\, D^{m}\widehat{Z}(q)=0,
\label{eq:TMLDEgam}
\end{equation}
which is again of regular singular type provided $g_0(c) \neq 0$. We can repeat the results of Section~\ref{sec:EVOAs} concerning TMLDE series solutions and the rationality of $c$ and $h$ noting that $Z_V(-1/\tau,c)$ (c.f. \eqref{eq:ZVS}) satisfies  \eqref{eq:TMLDEgam} for 
$\left(\begin{smallmatrix}
\alpha&\beta\\ \gamma&\delta
\end{smallmatrix}\right)=
\left(\begin{smallmatrix}
0 & 1\\ 
-1 & 0
\end{smallmatrix}\right)$.
We therefore find the VOSA analogues of Proposition~\ref{prop:Zsol} and 
\ref{prop:Zgen}
\begin{proposition}
\label{prop:Zsoltwist}
Let $V$ be an Exceptional VOSA of lowest primary weight $l\in\N+\half$  and central charge $c$ and let $N$ be a simple ordinary $V$--module of lowest weight $h$. Assuming $g_{0}(c)\neq 0$ in the TMLDE  \eqref{eq:TMLDE} then
\begin{enumerate}
\item[(i)]   $Z_{N}(q)$ is  holomorphic in $q^{\half}$ for $0<|q|<1$. 
\item[(ii)]   
$Z_{N}\left(\frac{\alpha\tau +\beta}{\gamma\tau +\delta}\right)$ is a solution of the TMLDE \eqref{eq:TMLDEgam} for all
 $\bigl(\begin{smallmatrix}
\alpha&\beta\\ \gamma&\delta
\end{smallmatrix} \bigr)\in \SL(2,\Z)$. 
\item[(iii)] The central charge $c$  and the lowest conformal weight $h$ are rational. 
\end{enumerate}
\end{proposition}
\begin{proposition}
\label{prop:Zgentwist}
Let $V$ be an Exceptional VOSA of lowest primary weight $l\in\N+\half$  and central charge $c$. Assuming that $g_{0}(c)\neq 0$ and that $m\le l-\half$ for any indicial root of the form $x=m-c/24$. We then find 
\begin{enumerate}
\item[(i)] $ Z_V(q)$ is the unique $q^{\half}$--series solution of the TMLDE  with leading $q$--expansion $Z_V(q)=Z_{V_\omega}(q)+O\left (q^{l-c/24}\right)$.
\item[(ii)] $\dim V_n$ is a rational function of $c$ for each $n\in\N+\half $.
\item[(iii)]  $V$  is generated by the space of  lowest weight primary vectors $\Pi_l$.
\end{enumerate}
\end{proposition}
We verify below for $l\le 17/2$ that $g_{0}(c)\neq 0$ and  that $m\le l-\half$ for any indicial $x=m-c/24$. 
We conjecture  these conditions hold in  general.

We can  construct two infinite series of  $p_l=1$ Exceptional VOSAs which we conjecture are examples.
\begin{proposition}
\label{prop:ADVOAtwist}
For each Virasoro minimal model with $h_{1,p-1}\in\N+\half$ there exists an Exceptional VOSA with one odd parity primary vector of  lowest  weight $l=h_{1,p-1}$ of $AD$-type
\begin{equation}
V=L\left(c_{p,q},0\right)\oplus L\left(c_{p,q},h_{1,p-1}\right).
\label{eq:Vsimcur}
\end{equation}
\end{proposition}
\begin{proposition}
\label{prop:W3kVOAtwist}
For each $k\in\N+\half$ for $k\ge \frac{3}{2}$ there exists an Exceptional VOSA ${\mathcal W}(3k)$ with one odd parity primary vector of  lowest  weight $3k$ and central charge $c_k=1-24k$.
\end{proposition}
Finally, similarly to Section~\ref{sect:genus 0}, with $G=\Aut(V)$ we have
\begin{proposition}
\label{prop:irredtwist}
Let $V$ be an Exceptional VOSA of class $\mathcal{S}^{2l+1}$ with primaries ${\Pi_l}$ of  lowest weight $l\in\N+\half$. If $\Pi_l$ is a reducible $G$--module then it is either an irreducible $G$--module or the direct sum of two isomorphic irreducible $G$--modules.
\end{proposition} 

\section{Exceptional SVOAs with Lowest  Primary Weight with $l\in\N+\half$ for $l\le \frac{17}{2}$}
\label{sect: VOSAex}
We now consider   examples of  Exceptional VOSAs  of lowest primary weight $l\le \frac{17}{2}$. 
We  denote  by $E_n=E_n(q)$ the Eisenstein series and $F_n=E_{n}
\left[\begin{smallmatrix}
1 \\ 
-1
\end{smallmatrix}\right] (q )$ the twisted Eisenstein series of weight $n$ 
appearing in the order $l+\half$ TMLDE \eqref{eq:TMLDE}. 
For $l\le \frac{3}{2}$ we find all $c,h\in \Q$ whereas for $\frac{5}{2}\le l \le \frac{17}{2}$ we find all $c,h\in \Q$ for which $p_l=\dim \Pi_l\le 500000$
found by computer algebra techniques. We obtain many examples of known exceptional VOAs such as 
the free fermion VOSAs and the Baby Monster VOSA $\VB ^\natural
=\Com(V^\natural,\omega_{\frac{1}{2}})$, the commutant of  $V^\natural$ with respect to a Virasoro vector of central charge $\half$ \cite{Ho1}. 
Some other such commutant  theories also arise.
\medskip

\noindent $\boldsymbol{[l=\frac{1}{2}].}$
Propositions~\ref{prop:MLDEtwist}--\ref{prop:Zgentwist}
imply that $Z(q)$ satisfies the 1st  order TMLDE 
\begin{eqnarray*}
D Z+cF_{2}Z=0.
\end{eqnarray*}	
But $F_2(q)
=\frac{1}{24}+2\sum_{r\in\N+\half}\frac{rq^r}{1-q^r}$ so that
$
Z(q)=\left(\frac{\eta\left(\tau/2\right)}{\eta(\tau)}\right)^{2c}
$
with $p_{1/2}=2c$. 
An Exceptional VOSA exists for all $ p_{1/2}=m\in \N$
given by the tensor product of $m$ copies of the free fermion VOSA $V(H,\Z+\half )\cong L\left(c_{3,4},0\right)\oplus L\left(c_{3,4},\frac{1}{2}\right)$.

\medskip

\noindent $\boldsymbol{[l=\frac{3}{2}].}$
$Z(q)$ satisfies a 2nd  order TMLDE  
\begin{eqnarray*}
\,D^2 Z+\frac{2}{17}F_{2} ( 5 c+22 ) \,DZ
+\frac{1}{34}c \left( 4(5c+22) F_{4}+17 E_{{4}}\right)\, Z=0,
\end{eqnarray*}
with indicial roots $x_1=-c/24, x_2=(7c+24)/408$ with iterative solution
\begin{eqnarray*}
Z_V(q)&=&q^{-c/24}(1-p_{3/2}{q}^{3/2}+(1+p_2)q^2-(p_{3/2}+p_{5/2}){q}^{5/2}\ldots),
\\
p_{3/2}&=&\frac {8c ( 5c + 22 )} {3(2 c-49)},\quad 
p_2={\frac { ( 5c + 22 )  ( 4 c+21 )  ( 10 c-7 ) }{ 2( c-33 )  ( 2 c-49 ) }},\\
p_{5/2}&=&-{\frac {136 c ( 5c + 22 )  ( 4 c+21
 )  ( 10 c-7 ) }{15 ( 2 c-83 )  ( c-33 )  ( 2 c-49 ) }}.
\end{eqnarray*}
For $2c=-2\mod 17$, the indicial roots differ by an integer leading to denominator zeros for $p_n$. 
The $c,h\in \Q$ solutions with possible VOAs are
\begin{center}
\renewcommand{\arraystretch}{1.7}
\begin{tabular}{|c|c|c|c| c|c|}
\hline
$c$ & $p_{3/2}$ &$p_{2}$ &$p_{5/2}$ &
 VOSA &$h\in \Q$\\
\hline
$-{\frac {21}{4}}$ & 1 & 0 & 0 & 
$L\left(c_{3,8},0\right)\oplus L\left(c_{3,8},\frac{3}{2}\right)$ & $0,-\frac{1}{4}$\\
\hline
${\frac {7}{10}}$ & 1 & 0 & 0 &
 $L\left(c_{4,5},0\right)\oplus L\left(c_{4,5},\frac{3}{2}\right)$ & $0,\frac{1}{10}$\\
\hline
$\frac{15}{2}$ & 35 & 119 & 238 &  
$\Com\left(V_{\sqrt{2}E_8}^{+}, \omega_{\frac{1}{2}}\right)$ & $0,\frac{1}{2}$ \\
\hline 
$16$ & 256 & 2295 & 13056 & 
 $V_{\BW_{16}}^{+}\oplus \left(V_{\BW_{16}}^{+}\right)_{3/2}$& $0,1$\\
\hline
${\frac {114}{5}}$ & 2432 & 48620 & 537472 & 
$\Com \left(\VB ^\natural, \omega_{\frac{7}{10}}\right)$& $0,\frac{7}{5}$ \\
\hline
${\frac {47}{2}}$ & 4371 & 96255 & 1139374 &
$\VB ^\natural$& $0,\frac {49}{34}*$\\
\hline
\end{tabular}
\end{center}
The $c=\frac{15}{2}=8-\half$ VOSA is the commutant of  $V_{\sqrt{2}E_8}^{+}$ with respect to a Virasoro vector of central charge $\half$ with $\Aut(V)=S_8(2)$ \cite{LSY} and 
$\VB ^\natural$ is the Baby Monster VOSA with $\Aut(\VB ^\natural)=\mathbb{B}$ \cite{Ho1}. In both cases, $p_{3/2}$ is odd and $\Pi_{3/2}$ is $\Aut(V)$--irreducible in agreement with Proposition~\ref{prop:irredtwist} \cite{Atlas}. 
The $c=\frac{15}{2}$ VOSA is the simple current extension of the Barnes-Wall Exceptional VOA by its $h=\frac{3}{2}$ module.
The $c=\frac{114}{5}=\frac{47}{2}-\frac{7}{10}$ VOSA is the commutant of  $\VB ^\natural$ with respect to a Virasoro vector of central charge $\frac{7}{10}$ \cite{HLY,Y}. In the later case, we expect  $\Aut(V)=2.^2E_6(2):2$, the maximal subgroup of $\mathbb{B}$, which has a 2432 dimensional irreducible representation \cite{Atlas}.
$\VB ^\natural$ is self-dual so that the $h=\frac{49}{34}$ TMLDE solution is not an irreducible character. 
\medskip

\noindent $\boldsymbol{[l=\frac{5}{2}].}$
$Z(q)$ satisfies a 3rd order TMLDE
\begin{eqnarray*}
 && ( 734 c+49 ) \,D^3 Z 
+27  ( 2 c-1 )  ( 7 c+68) F_2 \,D^2 Z 
\\
&&
+ \left( 6  ( 7 c + 68 )  ( 2 c-1 )  ( 5 c + 22) F_4
+\frac{1}{2}  ( 2634 c^2+1739 c-29348 )  E_4\right) \,DZ 
\\
&&
+\Big(2 c ( 7c + 68 )  ( 2 c-1 )  ( 5 c + 22 ) F_6 
+\frac {27}{2} c  ( 2 c-1 ) 
 ( 7c + 68 ) E_4 F_2  
\\
&&
+5 c( 36 {c}^{2}+622 c
-2413 )  E_6 \Big)\, Z 
=0,
\end{eqnarray*} 
where 
\begin{eqnarray*}
p_{5/2}&=&{\frac { 8\left( 7 c+68 \right)  \left( 2 c+5 \right)  \left( 2
 c-1 \right)  \left( 5 c +22\right) }{5(8 c^3-716 c^2+16102 c+239)}}.
\end{eqnarray*}
There is one $c,h\in \Q$ solution with possible VOSA for $p_{5/2}\le 500000$
\begin{center}
\renewcommand{\arraystretch}2
\begin{tabular}{|c|c |c|c|}
\hline
$c$ & $p_{5/2}$ & VOSA & $h\in \Q$\\
\hline
$-\frac{13}{14}$ & 1 & $L\left(c_{4,7},0\right)\oplus L\left(c_{4,7},\frac{5}{2}\right)$ &
$0,-\frac{1}{14},\frac{1}{7} $ \\
\hline
\end{tabular}
\end{center}

\medskip

\noindent $\boldsymbol{[l=\frac{7}{2}].}$
$Z_V(q)$ satisfies a 
4th order TMLDE 
where $p_{7/2}=\frac{r(c)}{ s(c)}$ for 
\begin{eqnarray*}
r(c)&=& 128{ ( 5 c +22   )   ( 3 c+46  )   ( 2 c-1  )   ( 14+c  )   ( 5 c+3  )   ( 7 c +68  ) },\\
s(c)&=& 7(160 c^5-31176 c^4+2015748 c^3-41830202 c^2\\
&&-92625711 c+1017681).
\end{eqnarray*}
The $c,h\in \Q$ solutions with possible VOSA for $p_{7/2}\le 500000$ are
\begin{center}
\renewcommand{\arraystretch}{1.7}
\begin{tabular}{|c|c |c|c|}
\hline
$c$ & $p_{7/2}$ & VOSA & $h\in \Q$\\
\hline
$-\frac{161}{8}$ & 1 & $L\left(c_{3,16},0\right)\oplus L\left(c_{3,16},\frac{7}{2}\right)$ & $0,-\frac{5}{8},-\frac{3}{4},-{\frac {7}{8}}$\\
\hline
$-\frac{19}{6}$ & 1 & $L\left(c_{4,9},0\right)\oplus L\left(c_{4,9},\frac{7}{2}\right)$ & $0,-\frac{1}{9},-\frac{1}{6},\frac{1}{6}$ \\
\hline

\end{tabular}
\end{center}

\medskip

\noindent $\boldsymbol{[l=\frac{9}{2}].}$
$Z_V(q)$ satisfies a 5th order TMLDE 
where $p_{9/2}=\frac{r(c)}{ s(c)}$ for 
\begin{eqnarray*}
r(c)&=& 
160  ( 3 c+46 )   ( 2 c-1 )   ( 5 c+3 )   ( 11 c+232 )   ( 68+7 c )   ( 40 {c}^{2}+1778 c+11025 ), \\
s(c)&=&
9(3200 {c}^{6}-1096320 {c}^{5}+140381096 {c}^{4}-7850716276 {c}^{3}
+149541921538 {c}^{2}\\
&&+829856821745 c+7484560125).
\end{eqnarray*}
The $c,h\in \Q$ solutions with possible VOSA for $p_{9/2}\le 500000$ are
\begin{center}
\renewcommand{\arraystretch}{1.7}
\begin{tabular}{|c|c|c|c|}
\hline
$c$ & $p_{9/2}$ & VOSA & $h\in \Q$\\
\hline
$-\frac{279}{10}$ & 1 & $L\left(c_{3,20},0\right)\oplus L\left(c_{3,20},\frac{9}{2}\right)$ & $0,-{\frac {7}{10}},-1,-{\frac {11}{10}},-\frac{6}{5}$\\
\hline
$-\frac{125}{22}$ & 1 & $L\left(c_{4,11},0\right)\oplus L\left(c_{4,11},\frac{9}{2}\right)$ & 
$0,-{\frac {3}{22}},-{\frac {5}{22}},-\frac{3}{11},\frac{2}{11}$ \\
\hline
$-\frac{7}{20}$ & 1 & $L\left(c_{5,8},0\right)\oplus L\left(c_{5,8},\frac{9}{2}\right)$ & 
$0,-\frac{1}{20},\frac{1}{4},{\frac {7}{10}},{\frac {891}{1850}}*$ \\
\hline
$-35$ & 1 & ${\mathcal W}(\frac{9}{2})$ & $0,-{\frac {11}{10}},-\frac{4}{3},-\frac{7}{5},-\frac{3}{2}$\\
\hline
\end{tabular}
\end{center}
The $c=-\frac{7}{20},h=\frac {891}{1850}$ TMLDE solution  is not an irreducible character.

\medskip

\noindent $\boldsymbol{[l=\frac{11}{2}].}$
$Z_V(q)$ satisfies a 
6th order TMLDE 
where $p_{11/2}=\frac{r(c)}{ s(c)}$ for 
\begin{eqnarray*}
r(c)&=&- 640 ( 13 c+350  ) ( 7 c+25  ) ( 11 c+232  ) ( 2 c-1  ) ( 3 c+46  )   ( 68+7 c  ) 
\\
&& . ( 5 c+3  ) ( 10
 c-7  )   ( 40  c ^ 2 +3586 c+50743  ) ,\\
s(c)&=&11(2240000  c ^ 9 -1185856000  c ^ 8 +249718385120  c ^ 7 -25848494429040  c ^ 6 \\
&&
+1266635173648176  c ^ 5 -18264666939042072  c ^ 4 -
336264778062263522  c ^ 3\\
&&
 -861021133326393167  c ^ 2 +
653498177653904632 c-9760778116675215).
\end{eqnarray*}
The $c,h\in \Q$ solutions with possible VOSA for $p_{11/2}\le 500000$ are
\begin{center}
\renewcommand{\arraystretch}{1.7}
\begin{tabular}{|c|c|c|c|}
\hline
$c$ & $p_{11/2}$ & VOSA & $h\in \Q$\\
\hline
$-\frac{217}{26}$ & 1 & $L\left(c_{4,13},0\right)\oplus L\left(c_{4,13},\frac{11}{2}\right)$
 &${0,-\frac{2}{13},-{\frac {7}{26}},-{\frac {9}{26}},-{\frac {5}{13}},\frac {5}{26}}$\\
\hline
\end{tabular}
\end{center}
\medskip

\noindent $\boldsymbol{[l=\frac{13}{2}].}$
$Z_V(q)$ satisfies a 
7th order TMLDE 
with $p_{13/2}=\frac{r(c)}{ s(c)}$ for 
\begin{eqnarray*}
r(c)&=&
4480 ( 13 c+350  )   ( 5 c+164 )   ( 7 c+25  )   ( 11 c+232  )   ( 3 c+46  )   ( 4 c+21  ) \\  
&& ( 5 c+3  )  ( 10 c-7  ) 
( 1120  c ^ 4 +187160  c ^ 3 +6889980  c ^ 2 +58079018 c-24165453  ) ,
\\
s(c)&=&
13(125440000  c ^ {11} -94806656000  c ^{10} +29650660755200  c ^ 9
 -4865828683343040  c ^ 8 \\
&&+
431531398085049664  c ^ 7 -18001596789986119984  c ^ 6 +
107049283968364390448  c ^ 5
\\
&&
 +9359034900957509468076  c ^ 4 +
76817948684836018331724  c ^ 3 
\\
&&
+155170276090966927173843  c ^ 2 -
81951451902336562695126 c
\\
&&
-7944030229978323194805).
\end{eqnarray*}
The $c,h\in \Q$ solutions with possible VOSA for $p_{13/2}\le 500000$ are
\begin{center}
\renewcommand{\arraystretch}{1.7}
\begin{tabular}{|c|c|c|c|}
\hline
$c$ & $p_{13/2}$ & VOSA & $h\in \Q$\\
\hline
$-\frac{611}{14}$ & 1 & $L\left(c_{3,28},0\right)\oplus L\left(c_{3,28},\frac{13}{2}\right)$
 & $0,-{\frac {11}{14}},-{\frac {19}{14}},-\frac{3}{2},$\\
& & & $-{\frac {12}{7}},-{\frac {25}{14}},-{\frac {13}{7}}$ \\
\hline
$-\frac{111}{10}$ & 1 & $L\left(c_{4,15},0\right)\oplus L\left(c_{4,15},\frac{13}{2}\right)$ & $0,-\frac{1}{6},-\frac{3}{10},$\\
& & & $-\frac{2}{5},-\frac {7}{15},-\frac{1}{2},\frac{1}{5}$ \\
\hline
\end{tabular}
\end{center}

\medskip

\noindent $\boldsymbol{[l=\frac{15}{2}].}$
$Z_V(q)$ satisfies an 8th order TMLDE 
where $p_{15/2}=\frac{r(c)}{ s(c)}$ for 
\begin{eqnarray*}
r(c)&=& -28672 ( 13 c+350 )  ( 5 c+164 )  ( 7 c+25 )  ( 11 c+232 ) 
\\
&& 
.  ( 3 c+46 )  ( 17 c+658 ) ( 4 c+21 )  ( 5 c+3 )  ( 10 c-7 ) 
\\
&&
. ( 560 c^4+146584 c^3+9082444 c^2+133381952 c-27346605 ) 
,
\\
s(c)&=&
21073920000 c^{12}
-21694120448000 c^{11}
+9524271218201600 c^{10}
\\
&&
-2298054501201632000 c^9
+325029065007052546624 c^8
\\
&&
-26081744761028079338944 c^7
+968808700001847281619664 c^6
\\
&&
+787299295625321246276560 c^5
-696312046814218010729784676 c^4
\\
&&
-7887852431045609558472152948 c^3
-21020840196255652876820528205c^2
\\
&&
+3455907491220404701398711750 c
+4568101033862110116156159375.
\end{eqnarray*}
The $c,h\in \Q$ solutions with possible VOSA for $p_{15/2}\le 500000$ are
\begin{center}
\renewcommand{\arraystretch}{1.7}
\begin{tabular}{|c|c|c|c|}
\hline
$c$ & $p_{15/2}$ & VOSA & $h\in \Q$\\
\hline
$-\frac{825}{16}$ & 1 & $L\left(c_{3,32},0\right)\oplus L\left(c_{3,32},\frac{15}{2}\right)$ & 
$0,-{\frac {13}{16}},-{\frac {23}{16}},-\frac{7}{4},$\\
& & & $-{\frac {15}{8}},-{\frac {
33}{16}},-{\frac {17}{8}},-{\frac {35}{16}}
$\\
\hline
$-\frac{473}{34}$ & 1 & $L\left(c_{4,17},0\right)\oplus L\left(c_{4,17},\frac{15}{2}\right)$ & 
$ 0,-{\frac {3}{17}},-{\frac {11}{34}},-{\frac {15}{34}},$\\
& & & $-{\frac {9}{17}},-{\frac {10}{17}},-{\frac {21}{34}}, {\frac {7}{34}}
$\\ 
\hline
$-\frac{39}{10}$ & 1 & $L\left(c_{5,12},0\right)\oplus L\left(c_{5,12},\frac{15}{2}\right)$ & $0,\frac{1}{2},{\frac {13}{10}},-\frac{1}{6},-\frac{1}{5},\frac{2}{15}$\\
\hline
$\frac{25}{28}$ & 1 & $L\left(c_{7,8},0\right)\oplus L\left(c_{7,8},\frac{15}{2}\right)$ & $0,\frac{1}{28},{\frac {3}{28}},{\frac {5}{14}},\frac{3}{4},{\frac {9}{7}}$\\
\hline
$-59$ & 1 & ${\mathcal W}(\frac{15}{2})$ & $0,-{\frac {13}{7}},-{\frac {21}{10}},-{\frac {31}{14}},$
\\
 & & & $-{\frac {12}{5}},-{\frac {17}{7}},-\frac{5}{2},-{\frac {67}{62}}*$\\
\hline
\end{tabular}
\end{center}
The $c=-59,h=-\frac {67}{62}$ TMLDE solution  is not an irreducible character \cite{F}.

\medskip

\noindent $\boldsymbol{[l=\frac{17}{2}].}$
$Z_V(q)$ satisfies a 9th order TMLDE. 
The only $c,h\in \Q$ solution with possible VOSA for $p_{17/2}\le 500000$ is
\begin{center}
\renewcommand{\arraystretch}{1.7}
\begin{tabular}{|c|c|c|c|}
\hline
$c$ & $p_{17/2}$ & VOSA & $h\in \Q$\\
\hline
$-{\frac {637}{38}}$ & 1 & $L\left(c_{4,19},0\right)\oplus L\left(c_{4,19},\frac{17}{2}\right)$ & 
$0,-{\frac {7}{38}},-{\frac {13}{38}},-{\frac {9}{19}}
,-{\frac {11}{19}},
$\\
& & & $-{\frac {25}{38}},-{\frac {27}{38}},-{\frac {14}{19
}},{\frac {4}{19}}$\\
\hline
\end{tabular}
\end{center}

\end{document}